\RequirePackage{fix-cm}

\documentclass{article}

\usepackage[T1]{fontenc}
\usepackage[utf8]{inputenc}
\usepackage{graphicx}
\usepackage{amsmath,amsfonts,amssymb,amsbsy,mathabx,amsthm}
\usepackage[svgnames]{xcolor}
\usepackage{tikz}
\usepackage{lmodern}
\usepackage{pgfplots}
\pgfplotsset{compat=newest}

\usepackage{geometry}
\DeclareFontShape{U}{matha}{m}{n}{<-6> matha5<6-7> matha6<7-8> matha7<8-9> matha8<9-10> matha9<10-12> matha10 <12-> matha12}{}
\DeclareFontShape{U}{mathb}{m}{n}{<-6> mathb5<6-7> mathb6<7-8> mathb7<8-9> mathb8<9-10> mathb9<10-12> mathb10 <12-> mathb12}{}
\DeclareFontShape{U}{mathx}{m}{n}{<-6> mathx5<6-7> mathx6<7-8> mathx7<8-9> mathx8<9-10> mathx9<10-12> mathx10 <12-> mathx12}{}
 
\usepackage{subcaption}
\usepackage{import}

\makeatletter
\newif\if@restonecol
\makeatother

\usepackage[algo2e,ruled,lined,titlenumbered,commentsnumbered]{algorithm2e}

\graphicspath{{figures/}}

\theoremstyle{remark}
\newtheorem*{rem}{Remark}

\newcommand{\un}[1]{\ensuremath{\underline{#1}}}

\newcommand{\R}{\ensuremath{\mathbb{R}}}

\newcommand{\eft}{\ensuremath{\mathbf{t}}}

\newcommand\shapef{\varphi}
\newcommand{\stiff}{\ensuremath{\mathbf{K}}}
\newcommand{\force}{\ensuremath{\mathbf{f}}}
\newcommand{\domain}{\ensuremath{\Omega}}

\newcommand{\res}{\ensuremath{\mathbf{r}}}
\newcommand{\bz}{\ensuremath{\mathbf{z}}}
\newcommand{\bw}{\ensuremath{\mathbf{w}}}
\newcommand{\bp}{\ensuremath{\mathbf{q}}}

\newcommand{\tlam}{\ensuremath{\boldsymbol{\delta{\lambda}}}}

\newcommand{\tdep}{\ensuremath{\mathbf{\delta{u}}}}

\newcommand\hooke{\mathbb{H}}

\newcommand{\zero}{\ensuremath{{\mathbf{0}}}}

\newcommand\norm[2]{\left\|#1\right\|_{#2}}
\newcommand\Ker[1]{\operatorname{Ker}\left(#1\right)}
\newcommand\Range[1]{\operatorname{Range}\left(#1\right)}
\newcommand\Rank[1]{\operatorname{Rank}\left(#1\right)}

\newcommand\divg{\ensuremath{\operatorname{div}}}

\newcommand{\Ltwo}{L^2}
\newcommand{\Honed}[1]{\left[H^1(#1)\right]^d}
\newcommand{\KA}{\mathcal{K}\!\mathit{a}}

\newcommand{\SA}{\mathcal{S}\!\mathit{a}}

\newcommand{\dom}{\Omega}
\newcommand{\du}{\partial_u}

\newcommand{\df}{\partial_f}
\newcommand{\dfdom}{\partial_f\Omega}

\newcommand\eps{\ensuremath{\varepsilon}}
\newcommand\sig{\ensuremath{\sigma}}
\newcommand\fvol{f}
\newcommand\fsurf{g}
\newcommand\young{\mathrm{Y}}

\newcommand\shapev{\boldsymbol{\shapef}_h}
\newcommand{\dep}{\ensuremath{\mathbf{u}}}
\newcommand{\Dep}{\ensuremath{\mathbf{U}}}

\newcommand{\citep}[1]{\cite{#1}}
\newcommand\Nsd{\mathcal{N}\!_{\mathrm{sd}}}
\newcommand{\s}{\ensuremath{^{(s)}}}
\newcommand{\sT}{\ensuremath{^{(s)^T}}}
\newcommand{\sinv}{\ensuremath{^{(s)^{-1}}}}

\newcommand{\ddl}{\ensuremath{}}
\newcommand{\ddc}{\ensuremath{}}

\newcommand{\ddlT}{\ensuremath{^{T}}}

\newcommand{\interf}{\ensuremath{\Upsilon}}
\newcommand{\bounda}{\ensuremath{\Gamma}}

\newcommand\trace{\ensuremath{\operatorname{tr}}}
\newcommand{\traceh}{\ensuremath{\mathbf{t}}}

\newcommand{\lam}{\ensuremath{\boldsymbol{\lambda}}}

\newcommand{\lamb}{\ensuremath{\boldsymbol{\lambda}_b}}

\newcommand{\LamF}{\ensuremath{\boldsymbol{\Lambda}}_F}

\newcommand{\passem}{\ensuremath{\mathbf{A}}}
\newcommand{\dassem}{\ensuremath{\mathbf{B}}}

\newcommand{\dassemF}{\ensuremath{\dassem_\mathrm{F}}}
\newcommand{\psassem}{\ensuremath{\mathbf{\tilde{A}}}}
\newcommand{\dsassem}{\ensuremath{\mathbf{\tilde{B}}}}

\newcommand{\Lam}{\ensuremath{\boldsymbol{\Lambda}}}

\newcommand{\Pop}{\ensuremath{\mathbf{P}}}

\newcommand\hu{\hat{u}}

\newcommand\hsig{\hat{\sig}}
\newcommand\ecr[1]{\boldsymbol{e}_{\mathrm{CR},#1}}
\newcommand\unorm[2]{|\!|#1|\!|_{\hooke,#2}}
\newcommand\signorm[2]{|\!|#1|\!|_{\hooke^{-1}\!,#2}}
\newcommand\hF{\hat{F}}
\newcommand\Feq{\mathcal{F}\!_{\mathrm{eq}}}

\newcommand\Ecal{\mathcal{E}}
\newcommand\Ncal{\mathcal{N}}

\newcommand\bF{\mathbf{F}}
\newcommand\bhF{\hat{\mathbf{F}}}

\newcommand{\kerMP}{\ensuremath{\mathbf{R}^\circlearrowright}}
\newcommand{\kerSP}{\ensuremath{\mathbf{v}^\circlearrowright}}

\newcommand{\kerMPT}{\ensuremath{\mathbf{R}^{{\footnotesize\circlearrowright}^T}}}

\colorlet{blue}{RoyalBlue!65!black}
\colorlet{green}{Green!65!black}
\colorlet{red}{red!65!black}
\colorlet{yellow}{Goldenrod!65!black}
\colorlet{violet}{Fuchsia!65!black}
\colorlet{biblio}{Fuchsia!45!black}
\colorlet{bghaut}{blue!0}
\colorlet{bgbas}{blue!10}
\tikzstyle{every picture}+=[align=center, remember picture]
\tikzstyle{base}=[rectangle, text centered, inner sep=5pt,
align=center,line width=.8pt]
\tikzstyle{baseb}=[rectangle, rounded corners=2pt, text centered, inner sep=5pt,
align=center,line width=.8pt]
\tikzstyle{basec}=[rectangle, rounded corners=2pt, text centered, inner sep=5pt,
align=center,line width=.8pt]
\tikzstyle{arrow}=[->,>=latex,thick]
\tikzstyle{redb}=[baseb,draw=red!65!black, fill=red!15, text=red!65!black]
\tikzstyle{blueb}=[baseb,draw=RoyalBlue!65!black, fill=RoyalBlue!15, text=RoyalBlue!65!black]
\tikzstyle{greenb}=[baseb,draw=Green!65!black, fill=Green!15, text=Green!65!black]
\tikzstyle{yellowb}=[baseb,draw=Goldenrod!65!black, fill=Goldenrod!30, text=Goldenrod!65!black]
\tikzstyle{violetb}=[baseb,draw=Fuchsia!65!black, fill=Fuchsia!15, text=Fuchsia!65!black]
\tikzstyle{greyb}=[baseb,draw=black!50, fill=black!20, text=black]
\tikzstyle{greybbb}=[base,draw=black, fill=black!20, text=black]
\tikzstyle{redt}=[text=red!65!black]
\tikzstyle{bluet}=[text=RoyalBlue!65!black]
\tikzstyle{greent}=[text=Green!65!black]
\tikzstyle{yellowt}=[text=Goldenrod!75!black]
\tikzstyle{violett}=[text=Fuchsia!65!black]
\tikzstyle{bluea}=[arrow,draw=RoyalBlue!65!black, fill=RoyalBlue!65!black]
\tikzstyle{reda}=[arrow,draw=red!65!black, fill=red!65!black]
\tikzstyle{greena}=[arrow,draw=Green!65!black, fill=Green!65!black]
\tikzstyle{violeta}=[arrow,draw=Fuchsia!65!black, fill=Fuchsia!65!black]
\tikzstyle{greya}=[arrow,draw=black!75, fill=black!75]
\tikzstyle{bluen}=[circle, minimum size=0.1, draw=RoyalBlue!65!black, fill=RoyalBlue!65!black]
\tikzstyle{redn}=[circle, minimum size=0.1, draw=red!65!black, fill=red!65!black]
\tikzstyle{greenn}=[circle, minimum size=0.1, draw=Green!65!black, fill=Green!65!black]
\tikzstyle{violetn}=[circle, minimum size=0.1, draw=Fuchsia!65!black, fill=Fuchsia!65!black\tikzstyle{greyn}=[circle, minimum size=0.1, draw=black!50, fill=black!50]
\tikzstyle{blackn}=[circle, minimum size=0.1, draw=black, fill=black]
\begin{document}






\title{Improved recovery of admissible stress in domain decomposition methods --- application to heterogeneous structures and new error bounds for FETI-DP}

\author{A.~Parret-Fr\'eaud$^2$, V.~Rey$^1$, P.~Gosselet$^1$, C.~Rey$^2$ \\[5pt]
$^1$ LMT-Cachan / ENS-Cachan, CNRS, Universit\'e Paris Saclay
\\61, avenue du pr\'esident Wilson, 94235 Cachan, France,\\[5pt]
$^2$ Safran Tech, rue des Jeunes Bois \\Ch\^ateaufort CS 80112,78772 Magny les Hameaux, France}

\maketitle

\begin{abstract}
This paper investigates the question of the building of admissible stress field in a substructured context. More precisely we analyze the special role played by multiple points. This study leads to (1) an improved recovery of the stress field, (2) an opportunity to minimize the estimator in the case of heterogeneous structures (in the parallel and sequential case), (3) a procedure to build admissible fields for FETI-DP and BDDC methods leading to an error bound which separates the contributions of the solver and of the discretization.

\end{abstract}




\section{Introduction}
\label{sec:intro}

The strong mathematical properties of the finite element method \cite{bre2008} for the approximation of the solution 
of mechanical problems are unfortunately not sufficient to precisely guarantee a priori the quality of the computed 
fields.
A posteriori verification aims at providing a numerical estimate of the distance between the unknown exact solution 
and 
the calculation. Several types of error estimators exist.

Estimators based on the lack of regularity \cite{Zie87} of the stress field are often efficient and easy to implement 
but they may underestimate the error.
Strict estimators can be obtained at the price of the evaluation of constants that depend on the shape of the domain, 
which is not very practical. Estimators based on the error in constitutive equation \cite{Lad75,Lad04}, or on the 
equilibrated residuals \cite{Kel84,Kel89} (which in fact are equivalent approaches) give strict bounds without 
constant but they require the computation of a statically  admissible (equilibrated) stress field.

Statically admissible stress field can be obtained by separate dual analysis \cite{Kem09}, or by a post-processing of 
the finite element displacement.
Among these post-processors, we can cite the Element Equilibration Technique (EET) \cite{Lad83,Lad97} and its 
recent variant \cite{Lad10bis,Ple11, VREY.2013.1}, and flux-free techniques based on partition of unity 
\cite{Car00,Par09,Par09bis,Gal09,Mau09,Moi09}.
Note that whatever the approach, the numerical cost is never negligible.

In \cite{Par10}, it was 
shown that the methods to post-process balanced stress fields could be embedded within the framework of 
non-overlapping domain decomposition \cite{Gos06}. In particular, when using the balancing domain 
decomposition (BDD \cite{Man93}) or the finite element 
tearing and interconnecting (FETI \cite{Far94bis}), specific displacement and traction fields can be generated and as inputs of a parallel procedure for the computation of statically admissible fields. Those results are recalled in 
Section~\ref{sec:principles}.

The first aim of this article, which is a subject of Section~\ref{sec:ptsmulti},  is to investigate more deeply the 
proposed procedure and in particular to emphasize the role played by nodes shared by several subdomains.
These nodes, often referred to as multiple points or cross-points, require specific attention for the 
computation of statically admissible stress fields. In particular, in case of strong heterogeneities, ignoring the 
importance 
of multiple points may lead to a defective error estimation. To tackle this difficulty, we analyze 
the role of the 
multiple points in error estimation and take advantage of the optimization problem it triggers to better the 
reconstruction of the admissible stress field. 
In Section~\ref{sec:seqhetero}, we show how the 
classical procedure EET to build statically admissible fields can be improved in order to take into account strong 
heterogeneity. The result is a more accurate sequential error estimator. We highlight the strong similarity between 
the optimization at the multiple points and the optimization in this sequential estimator. Then in Section~\ref{sec:numass}, 
we give numerical results of the optimized procedures applied to a two-dimensional mechanical problem with strong heterogeneities. 
Finally, we extend the parallel error estimation procedure to the FETI-DP algorithm \cite{Far00bis, Far01} and BDDC 
algorithm by a 
specific treatment of multiple points for the construction of balanced 
tractions. The bound separates the contributions of the solver and of the discretization.



\section{Principle of error estimation in substructured problems}
\label{sec:principles}

This section recalls the main principles of our approach for a posteriori error estimation in the framework of non-overlapping domain decomposition, with application to FETI \cite{Far94bis} and BDD \cite{Man93} as presented in \cite{Par10}.

\subsection{Reference mechanical problem and  admissibility spaces}
\label{sec:refmechpb}

Let us consider the static equilibrium of a structure which occupies the open polyhedral domain $\dom\subset\mathbb{R}^d$ and which is submitted to given body force $f$, traction force $g$ on $\partial_f\dom$ and 
displacement $u_d$ on the complementary part $\partial_u\dom$ such that $\operatorname{meas}(\partial_u\dom)\neq0$.
We assume the structure undergoes small perturbations and that the material is linear elastic, characterized by the 
Hooke's tensor $\hooke$. In the following,  $u$ is the unknown displacement field, $\eps(u)$ is the symmetric part of 
the gradient and $\sigma$ is the Cauchy stress tensor.

For an open subset $\omega \subset \dom$, we note $\partial_u\omega=\partial_u\dom\cap\partial\omega$ and 
$\partial_f\omega=\partial_f\dom\cap\partial\omega$, and we introduce the following subspaces of    \emph{kinematically admissible} $(\KA)$ and  \emph{statically admissible} $(\SA)$ fields:
\begin{subequations}
  \begin{align}
    \label{eq:fvCA}
    \KA(\omega)&=\left\{v\in\Honed{\omega}\ \text{such that}\ \trace(v)_{|\du\omega}=u_d\right\},\\
        \label{eq:fvCA0}
    \KA^0(\omega)&=\left\{v\in\Honed{\omega}\ \text{such that}\ \trace(v)_{|\du\omega}=0\right\},\\
        \label{eq:fvCA00}
    \KA^{00}(\omega)&=\left\{v\in\Honed{\omega}\ \text{such that}\ 
\trace(v)_{|\partial\omega\setminus\df\omega}=0\right\},\\
    \label{eq:fvSA}
    \SA(\omega)&= \left\{\begin{aligned} & \tau\in\left[\Ltwo(\omega)\right]^{d\times d}_{sym}  \text{ such that}\ 
\forall\ v\in\KA^{00}(\omega)\\  &\int_{\omega}\tau:\eps(v)\,dx = \int_{\omega} \fvol.v\,dx + \int_{\partial_f\omega} 
\fsurf.v\,dS\end{aligned}\right\}
  \end{align}
\end{subequations}
where $\trace$ is the trace operator. 
We introduce the following functional called \emph{error in constitutive relation} 
for a pair of displacement and stress fields
$(v,\tau)\in\Honed{\omega}\times\left[\Ltwo(\omega)\right]^{d^2}_{sym}$:
\begin{align}
  \label{eq:Rdc}
  \ecr{\omega}(v,\tau)=\signorm{\tau-\hooke:\eps(v)}{\omega},\quad\text{with } 
\signorm{\bullet}{\omega}=\left[\int_{\omega}\bullet:\hooke^{-1}:\bullet\,dx\right]^{1/2}
\end{align}

The mechanical problem set on $\dom$ may be formulated as:
\begin{equation}
  \label{eq:pbRef}
  \text{Find}\ (u,\sig)\in\KA(\dom)\times\SA(\dom)\ \text{verifying}\  \ecr{\Omega}(u,\sig)=0,
\end{equation}
Under the given assumptions, the solution to this problem exists and is unique, so that $\sig$ matches the singleton defined by $\SA(\dom)\cap\hooke:\eps\left(\KA(\dom)\right)$.
\medskip

Let us now consider a decomposition of domain $\dom$ into open subsets $(\dom\s)_{1\leqslant s\leqslant \Nsd}$ ($\Nsd$ 
is the number of subdomains) so that $\dom\s\cap\dom^{(s')}=\emptyset$ for $s\neq s'$ and $\bar{\dom}=\cup_s 
\bar{\dom}\s$, and mark with superscript $(s)$ the restriction on domain $\dom\s$.
Under this framework, kinematic and static admissibility on the whole structure may be restricted to each 
sub-structure 
$\dom\s$ providing the verification of interface conditions, namely displacements continuity (\ref{eq:KAss}) and 
 balance of tractions (or action-reaction principle) (\ref{eq:SAss}). Therefore, we have, for any globally admissible pair 
$(v,\tau)$:
\begin{subequations}
  \begin{equation}
    \label{eq:KAss}
    v \in \KA(\dom) \Leftrightarrow 
    \left\lbrace \begin{array}{l} 
        v\s\in \KA(\dom\s),\ \forall s \\ 
        \trace(v\s)=\trace(v^{(s')}) \text{ on } \interf^{(s,s')},\ \forall(s,s')
      \end{array}\right.,
  \end{equation}
  \begin{equation}
    \label{eq:SAss} 
    \tau \in \SA(\dom) \Leftrightarrow 
    \left\lbrace \begin{array}{l} 
        \tau\s\in \SA(\dom\s),\ \forall s \\ 
        \tau\s\cdot n\s+\tau^{(s')}\cdot n^{(s')}=0 \text{ on }\interf^{(s,s')},\ \forall(s,s')
      \end{array}\right.,
  \end{equation}
\end{subequations}
where $\interf^{(s,s')}=\partial\dom\s\cap\partial\dom^{(s')}$ is the interface between $\dom\s$ and $\dom^{(s')}$.

Of course, the error in constitutive relation can be written as a sum of local contributions:
\begin{equation*}
\forall (v,\tau)\in\Honed{\omega}\times\left[\Ltwo(\omega)\right]^{d^2}_{sym}, 
  \left.\ecr{\dom}(v,\tau)\right.^2 = \sum_{s=1}^{\Nsd} \left(\ecr{\dom\s}(v\s,\tau\s)\right)^2.
\end{equation*}

\subsection{Basics on error estimation}
\label{sec:referror}

Our error estimation technique is based on the following property, called the Pragger-Synge theorem~\cite{Lad75}:
\begin{align}
  \label{eq:eqRdcBSup}
  \forall (\hat{u},\hat{\sigma})\in \KA(\dom)\times \SA(\dom),\
 \unorm{ \eps(u-\hat{u})}{\Omega}^2  + \signorm{\sig-\hat{\sig}}{\Omega}^2     = \ecr{\dom}^2(\hat{u},\hat{\sig})
\end{align}
where $(u,\sig)$ solves the reference problem \eqref{eq:pbRef}. We choose to measure the error in displacement $e=u-\hat{u}$, and introduce the following norm: 
$\vvvert e \vvvert_\dom:= \unorm{\eps(e)}{\Omega} \leqslant \ecr{\dom}(\hat{u},\hat{\sig})$.
Thus for any admissible approximation of the solution $(\hat{u},\hat{\sigma})\in \KA(\dom)\times \SA(\dom)$, the error in constitutive relation is a computable upper bound of the error. Then the problem of the error estimation can be addressed through the ability to build 
kinematically and statically admissible approximations.\medskip

In the case of the finite element approach (of the monolithic problem), the displacement approximation $u_h$ solves the 
following problem:
\begin{align}
  \text{Find}\ u_h\in\KA_h(\dom)&\ \text{such that}\  \sig_h=\hooke:\eps(u_h)\ \text{satisfies} \nonumber\\
  \label{eq:fvPbApp}
  \int_{\dom}\sigma_h:\eps(v_h)\,dx &= \int_{\dom}
  \fvol.v_h\,dx + \int_{\dfdom} \fsurf.v_h\,dS,\ \forall\
  v_h\in\KA_h^0(\dom).
\end{align}
where  $\KA_h(\dom)$ is a finite dimension subspace of $\KA(\dom)$ so that we set $\hat{u}=u_h$. Note that $u_h$ takes 
the following form  $u_h=\shapev\dep$ where $\shapev$ is the matrix of shape functions and  $\dep$ the vector of nodal 
displacements.

\medskip
The solution $\sig_h=\hooke:\eps(u_h)$ unfortunately does not belong to $\SA(\dom)$, except in trivial cases. The 
construction of the admissible stress field $\hsig_h\in\SA(\dom)$ requires to apply a whole procedure (often referred 
to as recovery of balanced residual). It is a crucial point since the sharpness of error estimation strongly depends on 
the quality of the $\hsig_h$. 

We note by $\Feq$ the algorithm chosen to build an admissible stress field $\hat{\sigma}_h$ from the finite element one 
$\sigma_h$. $\Feq$ also takes as inputs the continuous representation of the imposed body force~$f$ and of the traction force~$g$.
\begin{equation*}
  \hat{\sigma}_h=\Feq(\sigma_h, f, g) \in \SA(\dom)
\end{equation*}
In our applications, the chosen algorithm is either the Element Equilibration Technique \cite{Lad83} or the 
Flux-free technique \cite{Par09} with p-refinement for elements problems (three degrees higher polynomial basis 
\cite{Bab94}). The general principles of the EET technique are recalled in section \ref{sec:seqhetero} where an 
improvement is proposed in order to take into account the material heterogeneity.

\subsection{Admissible field recovery for FETI and BDD domain decomposition}
\label{sec:errorestim}

We assume that the mesh and the decomposition of $\dom$ are conforming so that (i) each element only belongs to one 
subdomain and (ii) nodes are matching on the interfaces. Under this assumption, each degree of freedom is either 
located inside a subdomain or on its boundary $\bounda\s=\cup_{s'}\interf^{(s,s')}$  where it is shared with at least 
one neighboring subdomain.

In order to decouple the subdomains, we introduce $\lamb\s$ the vector of unknown nodal reactions imposed on the 
boundary of $\dom\s$ by its neighbors. The finite element equilibrium of subdomain $s$ then writes:
\begin{equation}\label{eq:ddm:discr:eq}
\stiff\s\dep\s=\force\s + \traceh\sT \lamb\s
\end{equation}
where $\stiff\s$ is the stiffness matrix, $\force\s$ is the vector of generalized forces and $\traceh\s$ is the 
discrete trace operator (which extracts the boundary values from a vector defined on the whole subdomain) and thus 
$\traceh\sT$ is  the extension by zero operator.

To complete the system, we need the discrete counterparts of the interface equations: continuity of 
displacement~\eqref{eq:KAss} and balance of forces~\eqref{eq:SAss}. This is done by the introduction of assembling 
operators $\passem$ and $\dassem$ (they will play an important role in our analysis and will be thoroughly described in 
the next section):
\begin{subequations}
  \begin{align}
   \sum_s\passem\s\lam_b\s&=\zero \label{eq:passlam}\\
    \sum_s\dassem\s\traceh\s\dep\s&=\zero  \label{eq:dassu}
  \end{align}
  \end{subequations}
Note that when there are only two subdomains, and thus one simple interface, these equations can be  written as: 
$\dep_b^{(1)}-\dep_b^{(2)}=\mathbf{0}$ and $\lamb^{(1)}+\lamb^{(2)}=\mathbf{0}$.

In \cite{Par10}, it was proved that when using the most classical iterative solvers for system 
(\ref{eq:ddm:discr:eq},\ref{eq:passlam},\ref{eq:dassu}), namely FETI \cite{Far94bis} and BDD \cite{Man93}, it is 
possible to build at no extra cost the following fields, whatever the convergence state of the solver:
\begin{itemize}
\item $\left(\dep_D^{(s)}\right)$: vectors of displacements  continuous at the interface $(\sum_s\dassem\s\traceh\s\dep_D\s=\zero)$ which define a globally admissible field: 
$\hu_D=(\shapev\s\dep_D\s)_s\in\KA(\Omega)$.
\item $\left(\lam_N^{(s)}\right)$: vectors of nodal reaction which are balanced at the interface, so that 
$\sum_s\passem\s\lam_N\s=\zero$, and which make the local equilibrium \eqref{eq:ddm:discr:eq} well posed Neumann 
problems.
\item $\left(\dep_N^{(s)}\right)$: vectors of displacements associated with the resolution of \eqref{eq:ddm:discr:eq} 
with imposed reactions~$\left(\lam_N^{(s)}\right)$. Let $u_N=(\shapev\s\dep_N\s)_s$, the field 
$\sigma_N=\left(\hooke:\eps(u_N)\right)$ is a global stress field which satisfies the finite element equilibrium.
\end{itemize}
As a consequence, it was proposed,  still in \cite{Par10}, to use $\sigma_N$ as a starting point for the computation of 
a global statically admissible stress field. In order to use classical recovery algorithms, one simply had to give a 
continuous representation $g_\bounda\s\in L^2(\bounda\s)$ of the traction field on the boundary compatible with the 
nodal reactions $\left(\lam_N\s\right)$:
\begin{equation}\label{eq:lamtog}
\lam_N\sT = \int_{\bounda\s} g_\bounda\s\cdot{\shapev}_b\s dS
\end{equation}
Then one simply used the recovery algorithm in parallel:
\begin{equation*}
  \hat{\sigma}_N= \left(\Feq\left(\sigma_N\s, f\s, g\s, g_\bounda\s(\lam_N\s)\right)\right)_s 
\end{equation*}
leading to the following inequality \cite{Par10}:
 \begin{equation}\label{eq:inegagus}
\vvvert u-\hu_D \vvvert^2 \leqslant \sum_{s=1}^{\Nsd} \left(\ecr{\dom\s}(\hu\s_D,\hsig\s_N)\right)^2
\end{equation}

In fact, there are cases where the above approach is not satisfactory:
\begin{itemize}
\item There is a risk that in the vicinity of multiple points $g_\bounda\s(\lam_N\s)$ looses its balance ($g_{\bounda|\Upsilon^{(s,s')}}\s+g_{\bounda|\Upsilon^{(s,s')}}^{(s')}\neq 0$), because the contributions of the two (or more) neighbors are not correctly distinguished. This causes a theoretical loss (never met in practice) of the exactness of the error bounding.
\item In the presence of heterogeneity, the parallel estimator may become unduly large compared to the sequential one.
\end{itemize}
The first problem is addressed in section~\ref{sec:ptsmulti}, the second one in section~\ref{sec:seqhetero} and assessments are presented in section~\ref{sec:numass}. An important by-product of these developments, exposed in section~\ref{sec:fetidp}, is the extension of all previous methods and results to FETI-DP and BDD-C algorithms which are very popular domain decomposition methods \cite{Far01,Doh03}.

\begin{rem}
In \cite{VREY.2013.1.2}, it was proved that the following bounds with separated contributions could also be derived:
 \begin{equation}\label{eq:sepsource}
 \begin{aligned}
\vvvert u-\hu_D \vvvert &\leqslant \|\mathbf{r}\|_M +  \sqrt{\sum_{s=1}^{\Nsd} 
\left(\ecr{\dom\s}(u\s_N,\hsig\s_N)\right)^2} \\ 
\vvvert u-u_N \vvvert&\leqslant \|\mathbf{r}\|_M + \sqrt{  \sum_{s=1}^{\Nsd} 
\left(\ecr{\dom\s}(u\s_N,\hsig\s_N)\right)^2}\\ 
\end{aligned}
\end{equation}
where $\|\mathbf{r}\|_M$ is a well-chosen norm of the solver's residual computed at each iteration, whereas the second term is governed by the discretization error.
\end{rem}



\section{Recovery of admissible fields in the presence of multiple points}
\label{sec:ptsmulti}

In this section we present how multiple points need to be taken into account when recovering stress fields in parallel. To do so, we first show the impact of multiple points on assembling operators, then we present the improved recovery procedure and finally we show how an opportunity is left to optimize the estimator.


\subsection{Assembling operators}
\label{sec:ptsmulti:assem}

The assembling operators aim at connecting neighboring subdomains together. They rely on various descriptions of the 
interface between subdomains. An example is provided in figure~\ref{fig:assem} where operators are concatenated in a 
row $\passem\ddl = \begin{pmatrix}\passem^{(1)}&\ldots&\passem^{(\mathcal{N}_{sd})}\end{pmatrix}$. In this example, we use a crosspoint of multiplicity 4 (the node is shared by 4 subdomains) but all developments in this section remain valid whatever the 
multiplicity.

Since we assume conforming meshes and discretization at the interface, we can use a description of the substructures in 
terms of degrees of freedom. In the following expressions, the boundary operator $\partial$ gives the degrees of freedom 
associated with nodes on the boundary of the subdomain. We can define different sets of degrees of freedom on the 
interface. In all cases, degrees of freedom where Dirichlet conditions are imposed are excluded:
\begin{equation}
  \begin{aligned}
&    \interf^{(i,j)} = (\partial\domain^{(i)}\cap\partial\domain^{(j)})\setminus \partial_u \domain, \text{ interface between }\Omega^{(i)} \text{ and }\Omega^{(j)}\\
&    \bounda^{(i)} = \partial\domain^{(i)} \setminus \partial \domain = \bigcup_{j}\interf^{(i,j)}, \text{ boundary of }\domain^{(i)} \\
&    \interf_p = \bigcup_{i}\interf^{(i)},\text{ primal interface}
  \end{aligned}
\end{equation}

The first way to connect subdomains together is to locate their boundary $(\bounda^{(i)})$ in the global interface $\interf_p$ (also called primal interface). 
To do so we introduce the inclusion map of the local boundary into the global interface, which we call the  primal assembly operator:
\begin{equation}
  \begin{aligned}
    \passem\s :&\ \R^{\bounda\s} \longrightarrow \R^{\interf_p}.
  \end{aligned}
\end{equation}
$\passem\s$ is a full rank matrix with $\#\interf_p$ rows and $\#\bounda\s$ columns, with one $1$ per column and zeros elsewhere.
Figure \ref{fig:assem:passem} shows an example for operators $\passem\s$ in the simple case of four substructures with a primal interface constituted by five nodes.

The second type of connection is realized using any set of matrices $(\dassem\s)$ such that:
\begin{equation}
  \begin{aligned}
    \sum_s\dassem\s{\passem\s}^T &=0,\\
    \Rank{\dassem\s}=\#\bounda\s.
  \end{aligned}
\end{equation}
 $\dassem\s$ has $\#\bounda\s$ columns and all $(\dassem\s)_s$ have the same number of rows which shall 
be sufficiently large for each degree of freedom of $(\bounda\s)$ to be involved at least once.
The rows of $(\dassem\s)_s$ are called connecting relations and they are represented by the set $\interf_d$ (also called dual 
interface). 

Let $\R^{\bounda\ddc}$ represent the space of vectors defined on the boundary of the subdomains; we have 
the following fundamental property:
\begin{equation}\label{eq:spacesplit}
  \R^{\bounda\ddc}=\Range{\passem\ddlT}\overset{\perp}{\oplus}\Range{\dassem\ddlT}.
\end{equation}
In other words, a vector defined on the boundary of subdomains can be decomposed in a unique way as a the sum of a continuous (primal) contribution and a balanced (dual) contribution.\medskip

From a practical point of view, there exists many ways to construct the dual assembly operators $(\dassem\s)$. We only present the most classical assembly operator used by solvers (simply written $\dassem\s$) and a second assembly operator more relevant for the recovery of admissible fields (written $\dassemF\s$), we omit other possibilities \cite{Jus97}. 

The most classical dual assembly operator relies on the use of the following description of the dual interface:
  \begin{equation*}
    \begin{aligned}
      \interf_d =\interf_d^C := \left( \interf^{(i,j)} \right)_{i<j}
    \end{aligned}
  \end{equation*}
  and $\dassem\s$ is the signed inclusion map  $\R^{\bounda\s} \longrightarrow \R^{\interf_d^C}$.
  The signs ensure that two neighbors have opposite contributions.
  With that construction, which is very easy to implement, even subdomains sharing one single node create a connection.
  This leads to many redundant rows in $\dassem\ddl$, namely $m(m-1)/2$ for a $m$-multiplicity node while only $(m-1)$ relationships are necessary to ensure the connectivity. However, it is well known that redundancies do not affect the resolution procedure \citep{Far94bis}.
  The classical dual assembly operator is illustrated on figure \ref{fig:assem:dassem}.\medskip

 A second way to define the dual assembly operator is to only keep,  in the connectivity table interfaces with non-zero measure: these are faces in 3D (and edges in 2D). 
  \begin{equation*}
    \begin{aligned}
      \interf_{d} = \interf_d^{F} := \left\{ \interf^{(i,j)},\ i<j, \text{ with }\operatorname{meas}( \interf^{(i,j)})\neq 0\right\}
    \end{aligned}
  \end{equation*} 
 This dual assembly operator, written $\dassemF\ddl$, 
is illustrated on figure \ref{fig:assem:dassem_f}. Note that there is still one redundancy per multiple point. A specific basis of $\operatorname{ker}(\dassemF\ddlT)$, associated with a cyclic stress, is illustrated in figure~\ref{fig:huLamb}, we write it down $\kerMP$. Note that building $\kerMP$ is a rather simple operation which relies on an analysis of the mesh connectivity and of the decomposition.

\subsection{Scaled assembling operators}
The scaled assembling operators $(\psassem\s)$ and $(\dsassem\s)$ play a role in the preconditioning of the domain decomposition methods. They are such that $\sum\s \psassem\s\passem\sT = \mathbf{I}$ and $\sum\s \dassem\s\dsassem\sT\dassem^{(j)} = \dassem^{(j)}$. In other words, scaled assembling operators are pseudo-inverses of assembling operators: $\psassem=\passem^{T^+}$ and $\dsassem=\dassem^{T^+}$. Many scaling being possible, in order to characterized them, they are associated with optimization problems:
\begin{equation}
\begin{aligned}
&(\dep\s)_s \text{ given in }\R^\bounda\ddc,\ \Dep = \sum \psassem\s \dep\s \Leftrightarrow \Dep = \arg\min_{\R^{\interf_p}} \|\dep\s  - \passem\sT \Dep  \|_\mathcal{K} \\
&(\lam\s)_s \text{ given in }\R^\bounda\ddc,\ \Lam = \sum \dsassem\s \lam\s \Leftrightarrow \Lam = \arg\min_{\R^{\interf_d}} \|\lam\s  - \dassem\sT \Lam  \|_{\mathcal{K}^{-1}} 
\end{aligned}
\end{equation}
where $\mathcal{K}$ and $\mathcal{K}^{-1}$ stand for well-chosen norm of $\R^\bounda\ddc$. Note that because of the redundancies at multiple points, $\dsassem\s$ is not fully characterized by previous system, but anyhow the mechanically consistent piece of information  $\dassem\sT \Lam$ is unique (it only depends on $\mathcal{K}^{-1}$).

\subsection{Building a continuous representation of balanced interefforts}

In this subsection, we show how a nodal intereffort $\Lam_F\in\R^{\interf_d^F}$ allows a more precise recovery of interface fluxes than 
\eqref{eq:lamtog}. Indeed we can build the piecewise function $g_\bounda\s = (g_F^{(s,s')})_{s'}$ solution to:
\begin{equation}\label{eq:lamtogloc}
\left(\dassem_F\sT\Lam_{F}\right)^T_{|\interf^{(s,s')}} = \int_{\interf^{(s,s')}} g_F^{(s,s')}\cdot{\shapev}_b\s dS
\end{equation}
Typically one can develop $g_F^{(s,s')}$ on the finite element shape functions of $\interf^{(s,s')}$ as was done in \cite{Par10} on $\Gamma^{(s)}$.
Then we automatically have $g_F^{(s,s')}= - g_F^{(s',s)}\in L^2(\interf^{(s,s')})$, and the classical recovery can be used in parallel on the subdomains:
\begin{equation*}
  \hat{\sigma}_N= \left(\Feq\left(\sigma_N\s, f\s, g\s, \left(g_F^{(s,s')}(\dassem_F\sT\Lam_F)\right)_{s'}\right)\right)_s 
\end{equation*}

In the following subsection, we show how $\Lam_F\in\R^{\interf_d^F}$ can be obtained, starting from very general data, which are easily obtained in BDD and FETI \cite{Par10}, but also in FETI-DP approach as explained in section~\ref{sec:fetidp}.

\begin{figure}[ht]
  \centering
  \begin{subfigure}{1.\textwidth}
  \begin{tikzpicture}[scale=.42]
  \draw[base] (0,5) rectangle (3,8);
  \draw (1.5,6.5) node {\small $\dom^{(1)}$};
  {\tiny
    \draw (2.3,7.5) node {$1^{(1)}$};
    \draw (2.3,5.5) node {$2^{(1)}$};
    \draw (0.7,5.5) node {$3^{(1)}$};
  }
  \foreach \y in {5,8} {
    \foreach \x in {0,3} {
      \draw[blackn] (\x,\y) circle (0.1);
    }
  }
  \draw[base] (0,0) rectangle (3,3);
  \draw (1.5,1.5) node {\small $\dom^{(2)}$};
  {\tiny
    \draw (0.7,2.5) node {$1^{(2)}$};
    \draw (2.3,2.5) node {$2^{(2)}$};
    \draw (2.3,0.5) node {$3^{(2)}$};
  }
  \foreach \y in {0,3} {
    \foreach \x in {0,3} {
      \draw[blackn] (\x,\y) circle (0.1);
    }
  }
  \draw[base] (5,0) rectangle (8,3);
  \draw (6.5,1.5) node {\small $\dom^{(3)}$};
  {\tiny
    \draw (5.7,0.5) node {$1^{(3)}$};
    \draw (5.7,2.5) node {$2^{(3)}$};
    \draw (7.3,2.5) node {$3^{(3)}$};
  }
  \foreach \y in {0,3} {
    \foreach \x in {5,8} {
      \draw[blackn] (\x,\y) circle (0.1);
    }
  }
  \draw[base] (5,5) rectangle (8,8);
  \draw (6.5,6.5) node {\small $\dom^{(4)}$};
  {\tiny
    \draw (7.3,5.5) node {$1^{(4)}$};
    \draw (5.7,5.5) node {$2^{(4)}$};
    \draw (5.7,7.5) node {$3^{(4)}$};
  }
  \foreach \y in {5,8} {
    \foreach \x in {5,8} {
      \draw[blackn] (\x,\y) circle (0.1);
    }
  }

  \draw[greya,<->] (0,3.2) -- (0,4.8);
  \draw (0.3,4) node {\tiny $2$};    
  \draw[greya,<->] (8,3.2) -- (8,4.8);
  \draw (7.7,4) node {\tiny $4$};    
  \draw[greya,<->] (3.2,0) -- (4.8,0);
  \draw (4,0.35) node {\tiny $3$};    
  \draw[greya,<->] (3.2,8) -- (4.8,8);
  \draw (4,7.65) node {\tiny $1$};    
  \draw[greya] (4,4) -- (3.2,3.2);
  \draw[greya] (4,4) -- (3.2,4.8);
  \draw[greya] (4,4) -- (4.8,3.2);
  \draw[greya] (4,4) -- (4.8,4.8);
  \draw (4,4.5) node {\tiny $5$};    
  
  \draw (20,4) node {$
    \passem\ddl=\footnotesize\begin{pmatrix}
      \begin{pmatrix}
        1&0&0\\0&0&1\\0&0&0\\0&0&0\\0&1&0
      \end{pmatrix}
      \begin{pmatrix}
        0&0&0\\1&0&0\\0&0&1\\0&0&0\\0&1&0
      \end{pmatrix}
      \begin{pmatrix}
        0&0&0\\0&0&0\\1&0&0\\0&0&1\\0&1&0
      \end{pmatrix}
      \begin{pmatrix}
        0&0&1\\0&0&0\\0&0&0\\1&0&0\\0&1&0
      \end{pmatrix}
    \end{pmatrix}
    $};
\end{tikzpicture}\caption{Primal assembly operator $\passem\ddl$}
    \label{fig:assem:passem}
  \end{subfigure}
  
  \begin{subfigure}{1.\textwidth}
  \begin{tikzpicture}[scale=.41]
  \draw[base] (0,5) rectangle (3,8);
  \draw (1.5,6.5) node {\small $\dom^{(1)}$};
  {\tiny
    \draw (2.3,7.5) node {$1^{(1)}$};
    \draw (2.3,5.5) node {$2^{(1)}$};
    \draw (0.7,5.5) node {$3^{(1)}$};
  }
  \foreach \y in {5,8} {
    \foreach \x in {0,3} {
      \draw[blackn] (\x,\y) circle (0.1);
    }
  }
  \draw[base] (0,0) rectangle (3,3);
  \draw (1.5,1.5) node {\small $\dom^{(2)}$};
  {\tiny
    \draw (0.7,2.5) node {$1^{(2)}$};
    \draw (2.3,2.5) node {$2^{(2)}$};
    \draw (2.3,0.5) node {$3^{(2)}$};
  }
  \foreach \y in {0,3} {
    \foreach \x in {0,3} {
      \draw[blackn] (\x,\y) circle (0.1);
    }
  }
  \draw[base] (5,0) rectangle (8,3);
  \draw (6.5,1.5) node {\small $\dom^{(3)}$};
  {\tiny
    \draw (5.7,0.5) node {$1^{(3)}$};
    \draw (5.7,2.5) node {$2^{(3)}$};
    \draw (7.3,2.5) node {$3^{(3)}$};
  }
  \foreach \y in {0,3} {
    \foreach \x in {5,8} {
      \draw[blackn] (\x,\y) circle (0.1);
    }
  }
  \draw[base] (5,5) rectangle (8,8);
  \draw (6.5,6.5) node {\small $\dom^{(4)}$};
  {\tiny
    \draw (7.3,5.5) node {$1^{(4)}$};
    \draw (5.7,5.5) node {$2^{(4)}$};
    \draw (5.7,7.5) node {$3^{(4)}$};
  }
  \foreach \y in {5,8} {
    \foreach \x in {5,8} {
      \draw[blackn] (\x,\y) circle (0.1);
    }
  }

  \draw[greya,>->] (0,3.2) -- (0,4.8);
  \draw (0.3,4) node {\tiny $2$};    
  \draw[greya,>->] (3,3.2) -- (3,4.8);
  \draw (2.7,4) node {\tiny $6$};    
  \draw[greya,<-<] (5,3.2) -- (5,4.8);
  \draw (5.3,4) node {\tiny $8$};    
  \draw[greya,<-<] (8,3.2) -- (8,4.8);
  \draw (7.7,4) node {\tiny $4$};    
  \draw[greya,<-<] (3.2,0) -- (4.8,0);
  \draw (4,0.35) node {\tiny $3$};    
  \draw[greya,<-<] (3.2,3) -- (4.8,3);
  \draw (4,2.65) node {\tiny $7$};    
  \draw[greya,<-<]  (4.8,5)  -- (3.2,5);
  \draw (4,5.35) node {\tiny $5$};    
  \draw[greya,<-<] (3.2,8) -- (4.8,8);
  \draw (4,7.65) node {\tiny $1$};    
  \draw[greya,<-<] (3.2,3.2) -- (4.8,4.8);
  \draw (4.6,4.1) node {\tiny $10$};    
  \draw[greya,>->] (4.8,3.2) -- (3.2,4.8);
  \draw (3.8,4.6) node {\tiny $9$};
  
  \draw (21.5,4) node {$
    \dassem\ddl=\footnotesize  \begin{pmatrix}
      \begin{pmatrix}
        1&0&0\\0&0&1\\0&0&0\\0&0&0\\0&-1&0\\0&1&0\\0&0&0\\0&0&0\\0&1&0\\0&0&0
      \end{pmatrix}
      \begin{pmatrix}
        0&0&0\\-1&0&0\\0&0&1\\0&0&0\\0&0&0\\0&-1&0\\0&1&0\\0&0&0\\0&0&0\\0&1&0
      \end{pmatrix}
      \begin{pmatrix}
        0&0&0\\0&0&0\\-1&0&0\\0&0&1\\0&0&0\\0&0&0\\0&-1&0\\0&1&0\\0&-1&0\\0&0&0
      \end{pmatrix}
      \begin{pmatrix}
        0&0&-1\\0&0&0\\0&0&0\\-1&0&0\\0&1&0\\0&0&0\\0&0&0\\0&-1&0\\0&0&0\\0&-1&0
      \end{pmatrix}
    \end{pmatrix}
    $};
\end{tikzpicture}\caption{Classic dual assembly operator $\dassem\ddl$}
    \label{fig:assem:dassem}
  \end{subfigure}  
  
    \begin{subfigure}{1.\textwidth}
  \begin{tikzpicture}[scale=.41]
  \draw[base] (0,5) rectangle (3,8);
  \draw (1.5,6.5) node {\small $\dom^{(1)}$};
  {\tiny
    \draw (2.3,7.5) node {$1^{(1)}$};
    \draw (2.3,5.5) node {$2^{(1)}$};
    \draw (0.7,5.5) node {$3^{(1)}$};
  }
  \foreach \y in {5,8} {
    \foreach \x in {0,3} {
      \draw[blackn] (\x,\y) circle (0.1);
    }
  }
  \draw[base] (0,0) rectangle (3,3);
  \draw (1.5,1.5) node {\small $\dom^{(2)}$};
  {\tiny
    \draw (0.7,2.5) node {$1^{(2)}$};
    \draw (2.3,2.5) node {$2^{(2)}$};
    \draw (2.3,0.5) node {$3^{(2)}$};
  }
  \foreach \y in {0,3} {
    \foreach \x in {0,3} {
      \draw[blackn] (\x,\y) circle (0.1);
    }
  }
  \draw[base] (5,0) rectangle (8,3);
  \draw (6.5,1.5) node {\small $\dom^{(3)}$};
  {\tiny
    \draw (5.7,0.5) node {$1^{(3)}$};
    \draw (5.7,2.5) node {$2^{(3)}$};
    \draw (7.3,2.5) node {$3^{(3)}$};
  }
  \foreach \y in {0,3} {
    \foreach \x in {5,8} {
      \draw[blackn] (\x,\y) circle (0.1);
    }
  }
  \draw[base] (5,5) rectangle (8,8);
  \draw (6.5,6.5) node {\small $\dom^{(4)}$};
  {\tiny
    \draw (7.3,5.5) node {$1^{(4)}$};
    \draw (5.7,5.5) node {$2^{(4)}$};
    \draw (5.7,7.5) node {$3^{(4)}$};
  }
  \foreach \y in {5,8} {
    \foreach \x in {5,8} {
      \draw[blackn] (\x,\y) circle (0.1);
    }
  }

  \draw[greya,>->] (0,3.2) -- (0,4.8);
  \draw (0.3,4) node {\tiny $2$};    
  \draw[greya,>->] (3,3.2) -- (3,4.8);
  \draw (2.7,4) node {\tiny $6$};    
  \draw[greya,<-<] (5,3.2) -- (5,4.8);
  \draw (5.3,4) node {\tiny $8$};    
  \draw[greya,<-<] (8,3.2) -- (8,4.8);
  \draw (7.7,4) node {\tiny $4$};    
  \draw[greya,<-<] (3.2,0) -- (4.8,0);
  \draw (4,0.35) node {\tiny $3$};    
  \draw[greya,<-<] (3.2,3) -- (4.8,3);
  \draw (4,2.65) node {\tiny $7$};    
  \draw[greya,<-<]  (4.8,5)  -- (3.2,5);
  \draw (4,5.35) node {\tiny $5$};    
  \draw[greya,<-<] (3.2,8) -- (4.8,8);
  \draw (4,7.65) node {\tiny $1$};    
  
  \draw (21.5,4) node {$
    \dassem_{\mathrm{F}}\ddl=\footnotesize\begin{pmatrix}
      \begin{pmatrix}
        1&0&0\\0&0&1\\0&0&0\\0&0&0\\0&-1&0\\0&1&0\\0&0&0\\0&0&0
      \end{pmatrix}
      \begin{pmatrix}
        0&0&0\\-1&0&0\\0&0&1\\0&0&0\\0&0&0\\0&-1&0\\0&1&0\\0&0&0
      \end{pmatrix}
      \begin{pmatrix}
        0&0&0\\0&0&0\\-1&0&0\\0&0&1\\0&0&0\\0&0&0\\0&-1&0\\0&1&0
      \end{pmatrix}
      \begin{pmatrix}
        0&0&-1\\0&0&0\\0&0&0\\-1&0&0\\0&1&0\\0&0&0\\0&0&0\\0&-1&0
      \end{pmatrix}
    \end{pmatrix}
    $};
\end{tikzpicture}\caption{Face-only dual assembly operator $\dassemF\ddl$}
    \label{fig:assem:dassem_f}
  \end{subfigure}  
  \caption{Examples of primal and dual assembly operators}\label{fig:assem}
\end{figure}

\begin{figure}
 \begin{tikzpicture}[scale=.45]
  \draw[base] (0,5) rectangle (3,8);
  \draw (1.5,6.5) node {\small $\dom^{(1)}$};
  {\tiny
    \draw (2.3,7.5) node {$1^{(1)}$};
    \draw (2.3,5.5) node {$2^{(1)}$};
    \draw (0.7,5.5) node {$3^{(1)}$};
  }
  \foreach \y in {5,8} {
    \foreach \x in {0,3} {
      \draw[blackn] (\x,\y) circle (0.1);
    }
  }
  \draw[base] (0,0) rectangle (3,3);
  \draw (1.5,1.5) node {\small $\dom^{(2)}$};
  {\tiny
    \draw (0.7,2.5) node {$1^{(2)}$};
    \draw (2.3,2.5) node {$2^{(2)}$};
    \draw (2.3,0.5) node {$3^{(2)}$};
  }
  \foreach \y in {0,3} {
    \foreach \x in {0,3} {
      \draw[blackn] (\x,\y) circle (0.1);
    }
  }
  \draw[base] (5,0) rectangle (8,3);
  \draw (6.5,1.5) node {\small $\dom^{(3)}$};
  {\tiny
    \draw (5.7,0.5) node {$1^{(3)}$};
    \draw (5.7,2.5) node {$2^{(3)}$};
    \draw (7.3,2.5) node {$3^{(3)}$};
  }
  \foreach \y in {0,3} {
    \foreach \x in {5,8} {
      \draw[blackn] (\x,\y) circle (0.1);
    }
  }
  \draw[base] (5,5) rectangle (8,8);
  \draw (6.5,6.5) node {\small $\dom^{(4)}$};
  {\tiny
    \draw (7.3,5.5) node {$1^{(4)}$};
    \draw (5.7,5.5) node {$2^{(4)}$};
    \draw (5.7,7.5) node {$3^{(4)}$};
  }
  \foreach \y in {5,8} {
    \foreach \x in {5,8} {
      \draw[blackn] (\x,\y) circle (0.1);
    }
  }

  \draw[greya,>->] (3,3.2) -- (3,4.8);
  \draw (2.7,4) node {\tiny $6$};    
  \draw[greya,<-<] (5,3.2) -- (5,4.8);
  \draw (5.3,4) node {\tiny $8$};    
  \draw[greya,<-<] (3.2,3) -- (4.8,3);
  \draw (4,2.65) node {\tiny $7$};    
  \draw[greya,<-<]  (4.8,5) -- (3.2,5) ;
  \draw (4,5.35) node {\tiny $5$};    

  \draw (1.8,4) node {\footnotesize $\kerMP$};    
  
  \draw (18,4) node {$
    \kerMP=\begin{pmatrix}
      0\\ 0\\ 0\\ 0\\ 1\\1\\1\\1
      \end{pmatrix}
    $};
\end{tikzpicture}\caption{Illustration of $ \kerMP$ in the framework of figure \ref{fig:assem:dassem_f} } \label{fig:huLamb}
 \end{figure}

\subsection{General methodology to compute nodal intereffort $\Lam_F$}
\label{sec:ptsmulti:trac}

Using nodal interefforts like $\Lam_F$ in $\operatorname{Range}(\dassem_F\ddl)$ makes it simple to recover admissible stress fields. Unfortunately, it does not corresponds to the classical representation of the interface: in the primal approach, operator $\passem\ddl$ is used; whereas in the dual approach operator $\dassem\ddl$ is used, which differs from $\dassem_F\ddl$ as soon as there exist nodes shared by more than 3 subdomains (in which case $\dassem\ddl$ features interactions between subdomains with a zero-measure interface).


Our starting point for the construction of $\Lam_F$  is a vector of nodal tractions $\lam_N\s$ defined on the boundary of substructures such that:
\begin{align*}
\sum_s \passem\s\lam_N\s = 0
\end{align*}
which can always be obtained  during FETI or BDD iterations \citep{Par10}: in FETI, we directly have $\lam_N\s=\dassem\sT \Lam$ where $\Lam$ is the main unknown of the problem; whereas in BDD, we have  $\lam_N\s=\lam\s-\psassem\sT\sum_j\passem^{(j)} \lam^{(j)}$ where $(\lam\s)_s$ are  the nodal reactions resulting from local Dirichlet problems. The construction of $(\lam_N\s)$ during FETI-DP iterations is discussed in section~\ref{sec:fetidp}.
%

We thus wish to compute an interaction vector $\Lam_F$ which corresponds to that distribution of effort:
\begin{equation}\label{eq:lamFtolamN}
\lam_N\s =  \dassem_F\sT\Lam_F
\end{equation}
Thanks to \eqref{eq:spacesplit}, this problem possesses solutions. Because of the redundancy of $\dassemF\ddl$ at 
multiple points, solutions are defined up to a member of the null space of $\dassemF\ddl$. The 
solution is made unique by imposing that $\|\LamF - \LamF^1\|_\mathcal{P}$  shall be minimal in $\R^{\interf^F_d}$ for 
a chosen norm defined by a symmetric positive definite matrix $\Pop$ and a reference nodal interaction $\LamF^1$. 

The  $\mathcal{P}$-norm defines a weighted version of $(\dassemF\s)$ which we write down 
 $(\widehat{\dassem}_F\s)$:
\begin{align}
  \label{eq:hdepb_sol}
  \Lam_F & =  \sum_s\widehat{\dassem}_F\s \lam_N\s + \left( \mathbf{I}-\sum_s\widehat{\dassem}_F\s \dassem_F\sT \right)\Lam_F^1 \\
\text{with }  \widehat{\dassem}_F\s & = \Pop^{-1} \dassem_F\s \left(\sum_j\dassemF^{(j)^T}\Pop^{-1}\dassemF^{(j)}\right)^+
\end{align}
$(\widehat{\dassem}_F\s)$ is a sparse matrix with dense blocks for each multiple point. It is often computationally more interesting 
and numerically more stable to seek $\Lam_F$ as a corrected guess: let $\Lam_F^0$ be any vector satisfying 
\eqref{eq:lamFtolamN} (typically obtained with $\Pop=\mathbf{I}$), we have:
\begin{equation}\label{eq:roundGamma}
\Lam_F = \Lam_F^0 - \kerMP(\kerMPT\Pop\kerMP)^{-1}\kerMPT\Pop\left(\Lam_F^0- \LamF^1\right)
\end{equation}


The choices of the  $\mathcal{P}$-norm and of the reference interaction $\Lam_F^1$ are crucial to obtain a good estimation of the error, in particular in the presence of heterogeneities (see Section~\ref{sec:numass} for assessments). From our experiments, we draw two conclusions: first the heterogeneity must be taken into account, second $\Lam_F^1$ must be chosen in agreement with the real mechanical state of the structure which can be estimated through the finite element Cauchy stress on the element edges $\sig_h\s\cdot n^{(s,s')}$. It comes out that this issue has a direct equivalent in the sequential Element Equilibration Technique (EET) when optimizing the force fluxes around an internal node (see \cite{VREY.2013.1.2} for a discussion on the potential closed stress fluxes in finite element models). Note that, to the authors' knowledge, heterogeneities had never been considered in papers related to the EET.

\begin{figure}[ht]
\centering 
\includegraphics[width=1.\textwidth]{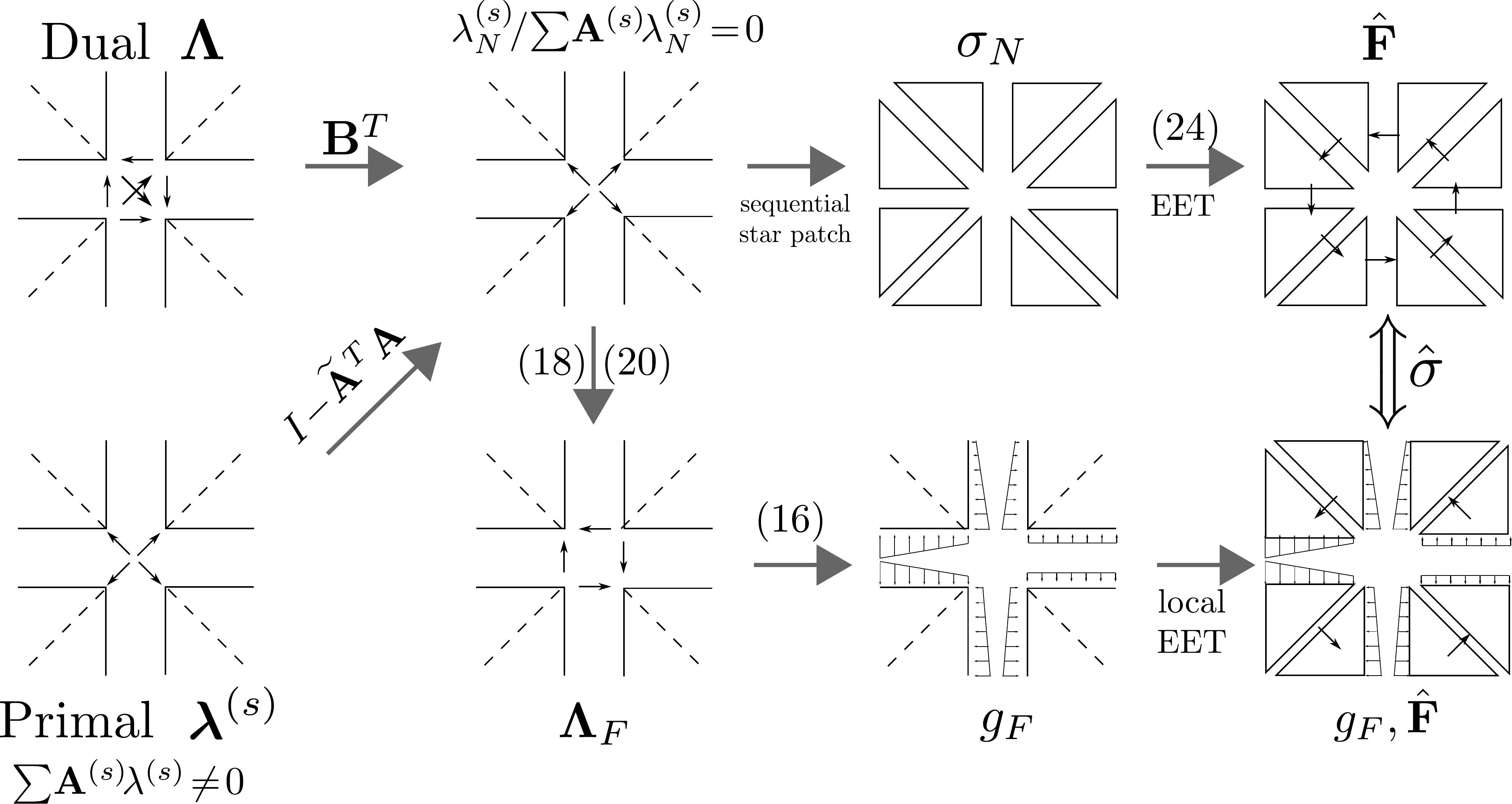}\caption{Methodology for parallel stress recovery / equivalent sequential star-patch}\label{fig:s1}
\end{figure}

In the next section, we present an improvement of the sequential EET in order to take into account heterogeneities. Regarding domain decomposition and multiple points, we propose to conduct the same computation. More precisely, we condense the star-patch problem built around multiple points on the interface in order to find the best $\Lam_F$. This strategy implies limited extra sparse communications: one small all-to-all exchange on local communicators associated with vertexes (in 2D) or edges (in 3D). The methodology to recover optimized stress is summed up in figure~\ref{fig:s1}.

%
%
%

\section{Improvement of sequential estimators for heterogeneous problems}
\label{sec:seqhetero}

This section presents how the Element Equilibration Technique (EET) can be improved to take into account heterogeneity. Moreover, as announced in previous section, these improvements have a direct counterpart in the parallel reconstruction of admissible stress field for domain decomposition methods. 

We first quickly recall the main principles of equilibrated stresses recovery through the EET algorithm and invite the interested reader to refer to \cite{Lad83,Lad93} for more details. Besides, it should be noted that the improvement presented here may also be applied to newer EESPT algorithm \cite{Lad10bis}, or the STARFLEET method \cite{VREY.2013.1}, since the common basis is shared by these approaches.\medskip

\subsection{Classical recovery of balanced stress}
The element-wise recovery of equilibrated stress takes place in two steps:
\begin{enumerate}
\item A traction field $(\hF_\gamma)_\gamma$ is built on the edges $\gamma$ (faces in 3D) of the mesh so that it verifies the equilibrium with both external and internal loading;
\item $\hsig_h$ is recovered from $(\hF_\gamma)$ through the resolution of element-wise problems defined as:
\begin{align*}
  \divg(\hsig_{h|E})+\fvol&=0\ &&\text{on }E, \\
  \hsig_{h|E}.n_E&=\delta_E^\gamma \hF_\gamma\ &&\text{on }\gamma\in\partial E,
\end{align*}
where $\delta^\gamma_E=\pm 1$ is a signed boolean coefficient verifying $\delta^\gamma_E+\delta^\gamma_{E'}=0$ on the 
edge $\gamma=\partial E\cap\partial E'$ in order to enable the balance of stress fluxes across element boundaries.
\end{enumerate}

As will be shown in section~\ref{sec:numass}, the jump of material coefficients need to be taken into account during the first step in order to build relevant values for equilibrated the 
traction fields $\hF_\gamma$, leading to effective error estimators.

The EET recovery method seeks the equilibrated traction field $(\hF_\gamma)$ thanks to the following prolongation condition on each element $E$ of the mesh:
\begin{align}
  \label{eq:eRdcSAEltProlStrong}
  \int_E(\hsig_h-\sig_h).\nabla\varphi_i\,dE=0,\ \forall E\in\Ecal_h,\ \forall
  i\in\Ncal_h^E,
\end{align}
where $\Ecal_h$ is the set of elements of $\dom_h$ and $\Ncal_h^E$ the set of vertexes of $E$.
This condition, which links the admissible stress field $\hsig_h$ to the finite element one $\sig_h$, leads to the following equations for $(\hF_\gamma)$:
\begin{equation}
  \label{eq:eRdcSAEltProlApplied}
  \sum_{\gamma\in\partial E}\int_{\partial 
E}\delta_E^\gamma \hF_\gamma\varphi_i\,d\Gamma=
\int_E\left(\sig_h:\eps(\varphi_i)-\fvol\varphi_i\right)dE,\nonumber \qquad 
  \forall E\in\Ecal_h,\ \forall i\in\Ncal_h^E,
\end{equation}

In a practical way, the previous system enables to seek for generalized nodal values $\bhF^\varphi_{\gamma,i}$ of 
$(\hF_\gamma)$ defined on each edge $\gamma$ adjacent to node $i$ by:
\begin{align*}
  \bhF^\varphi_{\gamma,i}=\int_{^\gamma}\delta_E^\gamma\hF_\gamma\varphi_id\Gamma.
\end{align*}

\begin{figure}[!ht]
  \centering
\includegraphics[width=.3\textwidth]{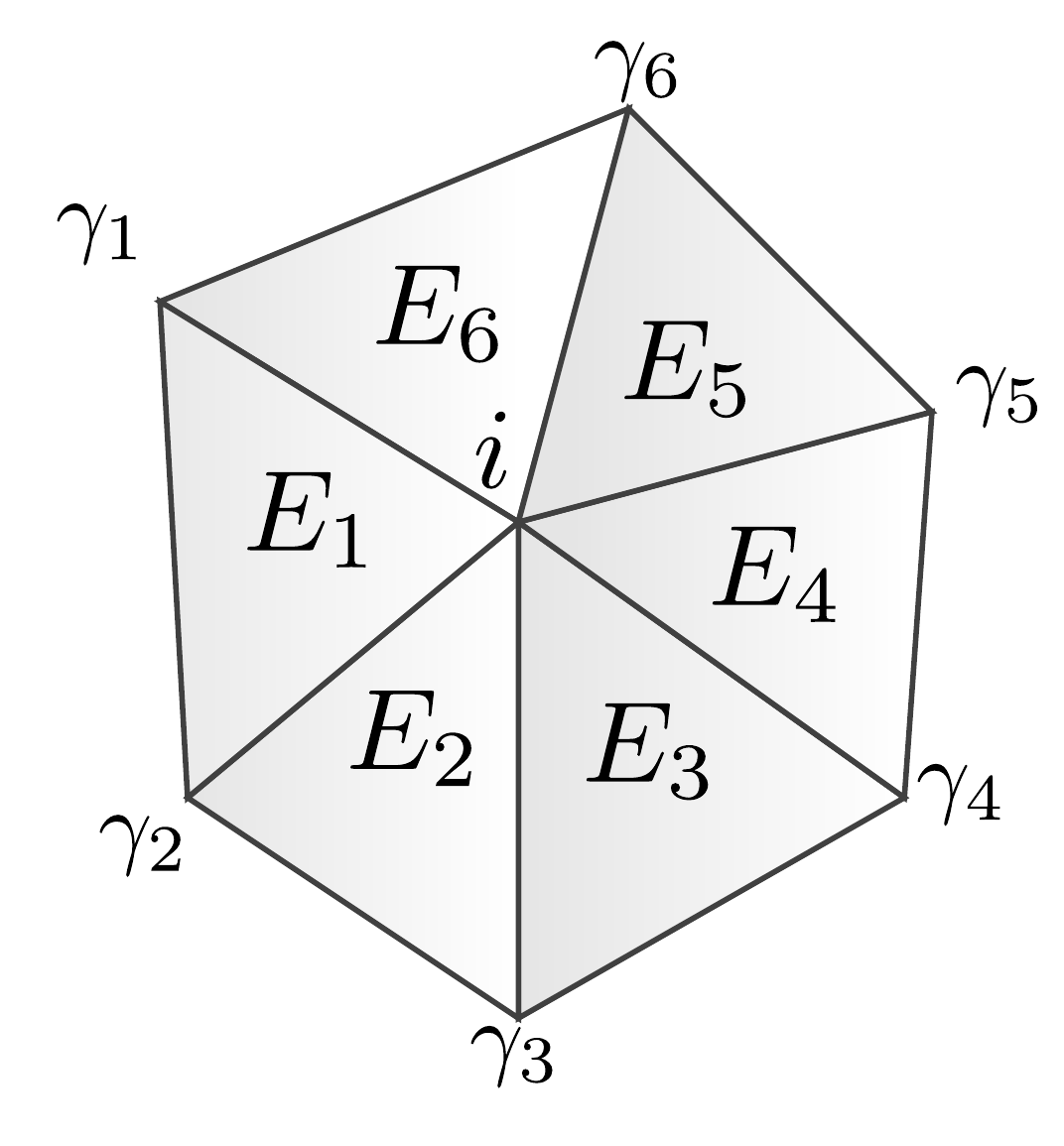}
  \caption{Star patch $\omega_{h,i}$ associated with the interior node $i$}
  \label{fig:EETpatch}
\end{figure}

We derive the resulting system in the case of one internal node as in figure~\ref{fig:EETpatch}. Note that similar optimization can be conducted in the case of nodes on the part of the boundary where a Dirichlet condition is prescribed (in other configurations the system is closed and leaves no opportunity for optimization). After expressing the condition \eqref{eq:eRdcSAEltProlApplied} for each element containing the node $i$, one obtains the following problem on patch $\omega_{h,i}$ corresponding to the support of shape function $\shapef_i$ (fig.~\ref{fig:EETpatch}):
\begin{equation}  \label{eq:EETsystPatch}
 \begin{pmatrix}
1 & -1 & 0  & 0  \\
0 & 1 & -1  & 0 \\
\vdots & \vdots & \ddots & \vdots\\
-1 & 0 & 0  & 1 
\end{pmatrix}
\begin{pmatrix}
  \bhF^\varphi_{\gamma_1,i} \\
\vdots \\    \bhF^\varphi_{\gamma_n,i}
\end{pmatrix}=
\begin{pmatrix}
\int_{E_1}(\sig_h:\eps(\varphi_i)-\fvol\varphi_i)dE \\ 
  \vdots \\ 
\int_{E_n}(\sig_h:\eps(\varphi_i)-\fvol\varphi_i)dE 
\end{pmatrix}
\end{equation}

The system \eqref{eq:EETsystPatch} is not sufficient to find an unique set of generalized nodal values $(\bhF^\varphi_{\gamma,i})_{\gamma\in\partial^-_i}$, where $\partial^-_i$ is the set of edges radiating from the node $i$. Indeed the vector $\kerSP=\begin{pmatrix}  1&\ldots&1\end{pmatrix}^T$ is a basis of the null space. It is clearly the sequential counterpart of the basis of $\Ker{\dassemF\ddlT}$, as illustrated on figure~\ref{fig:huLamb}.  
Classically, the solution is made unique by recasting the system in a constrained optimization problem:
\begin{equation}
\begin{aligned} 
& \text{minimize } \sum_{\gamma\in\partial^{-}_i}\left(\frac{\bhF^\varphi_{\gamma,i}-\bF^m_{\gamma,i}}{\operatorname{meas}(\gamma)}\right)^2 \text{under condition }\eqref{eq:EETsystPatch} \\
&\text{where }
\bF^m_{\gamma,i}=\int_{\gamma}\frac{1}{2}\left[\delta^\gamma_E\sig_{h|E}.n_E+\delta^\gamma_{E'}\sig_{h|E'}
.n_{E'} \right]
\varphi_i\,d\Gamma.
\end{aligned}
\end{equation}
From a mechanical point of view, previous problem consists in computing the  equilibrated traction field $(\bhF_{\gamma,i})_\gamma$ nearest to the mean values of stress fluxes across elements' edges  defined by $(\bF^m_{\gamma,i})_\gamma$.

\subsection{Improvement for heterogeneous problems}

In the heterogeneous case, considering for instance problems with strong gaps of the Young modulus across elements, the previous choice of minimization does not seem to be relevant. Following the strategies set up when preconditioning heterogeneous problems solved by FETI or BDD algorithm, we propose 
to use weighted mean stress with regard to Young modulus $Y$ in the minimization process:
\begin{equation}\label{eq:EETmoyHetero}
\begin{aligned} 
& \text{minimize } \sum_{\gamma\in\partial^{-}_i}p_{\gamma}\left(\frac{\bhF^\varphi_{\gamma,i}-\bF^m_{\gamma,i}}{\operatorname{meas}(\gamma)}\right)^2 \text{under condition }\eqref{eq:EETsystPatch} \\
&\text{where }
\tilde{\bF}^m_{\gamma,i}=\int_{\gamma}\frac{1}{\young_E^{-1}+\young_{E'}^{-1}}\left[\young_E^{-1}
\delta^\gamma_E\sig_{h|E}.n_E+\young_{E'}^{-1}\delta^\gamma_{E'}\sig_{h|E'}.n_{E'}\right]\varphi_i\,d\Gamma \\
&\text{and } p_{\gamma}=\left(\frac{1}{\young_{E'}}+\frac{1}{\young_E}\right)\  \text{ with }\gamma=\partial E\cap\partial E'.
\end{aligned}
\end{equation}
Therefore, the equilibrated traction field is closer to  the stress flux associated with the adjacent element with stronger flexibility.

\medskip
 Note that $(\bhF^\varphi_{\gamma,i})\gamma$ can be searched as a sum of a particular solution $(\bhF^{\varphi,P}_{\gamma,i})\gamma$ and of an element of the kernel $\beta\kerSP$ where $\beta$ is a scalar computed in order to minimize the following 
distance : 
\begin{align}
  \frac{1}{2}\norm{ \left(\bhF^{\varphi,P}_{\gamma,i} + \beta\kerSP -  \tilde{\bF}^m_{\gamma,i} \right)_{\gamma\in\partial^-_i}}{\mathcal{P}}
\end{align}
where $\mathcal{P}$ is the norm defined by $\operatorname{diag}(p_\gamma/\sqrt{\operatorname{meas}(\gamma)})$ which makes this minimization problem equivalent to~\eqref{eq:EETmoyHetero}. One can recognize the problem solved in the case of a multiple point~\eqref{eq:roundGamma} as detailed in the previous section.




\section{Numerical assessments}
\label{sec:numass}

\begin{figure}[ht]
  \centering
 \hfill \begin{subfigure}{.25\textwidth}
  \includegraphics[width=1.\textwidth]{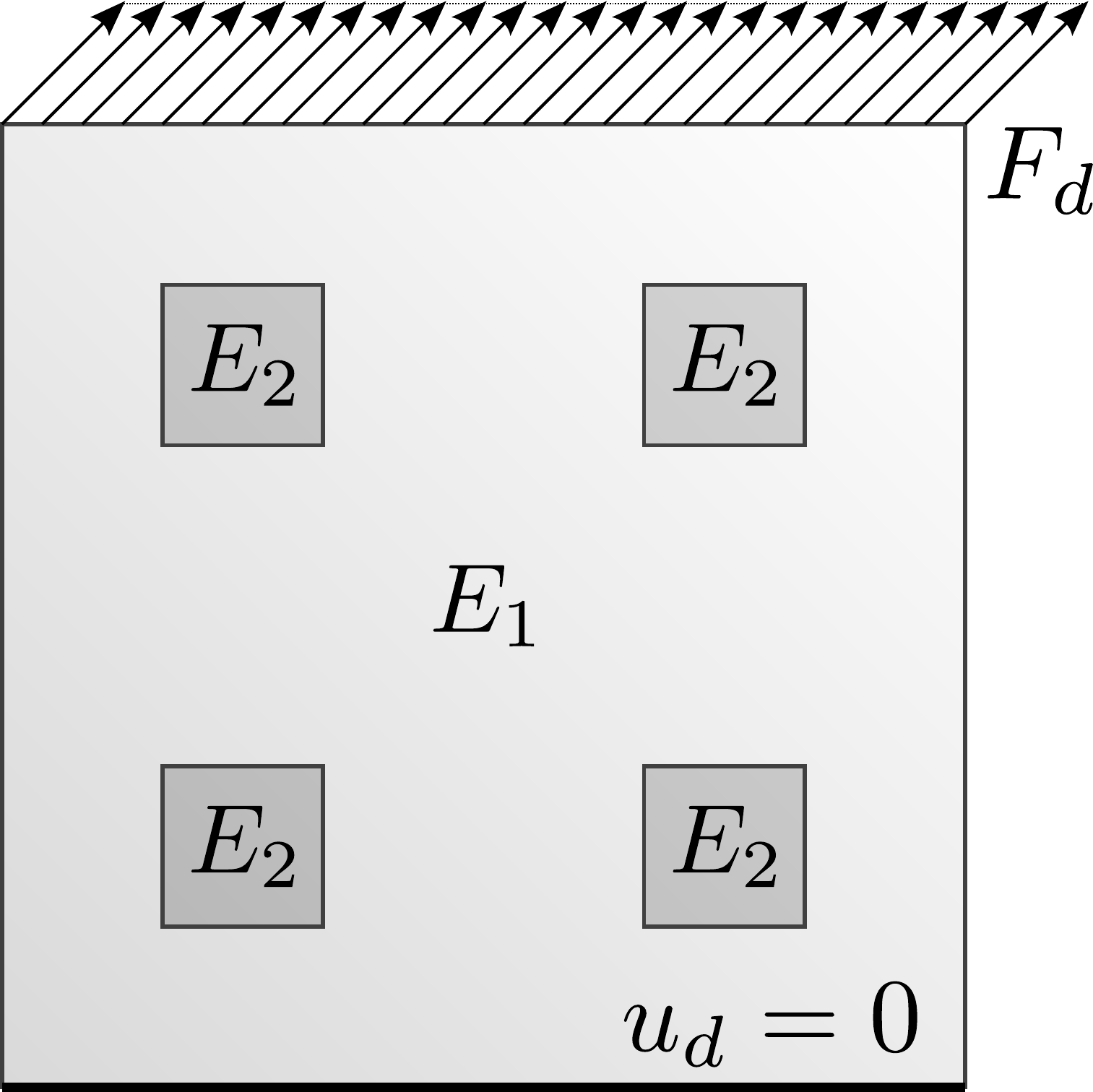}\caption{Description and loading}\label{fig:heteroPb}\hfill
\end{subfigure}\hfill
\begin{subfigure}{.24\textwidth}
\includegraphics[width=1.\textwidth]{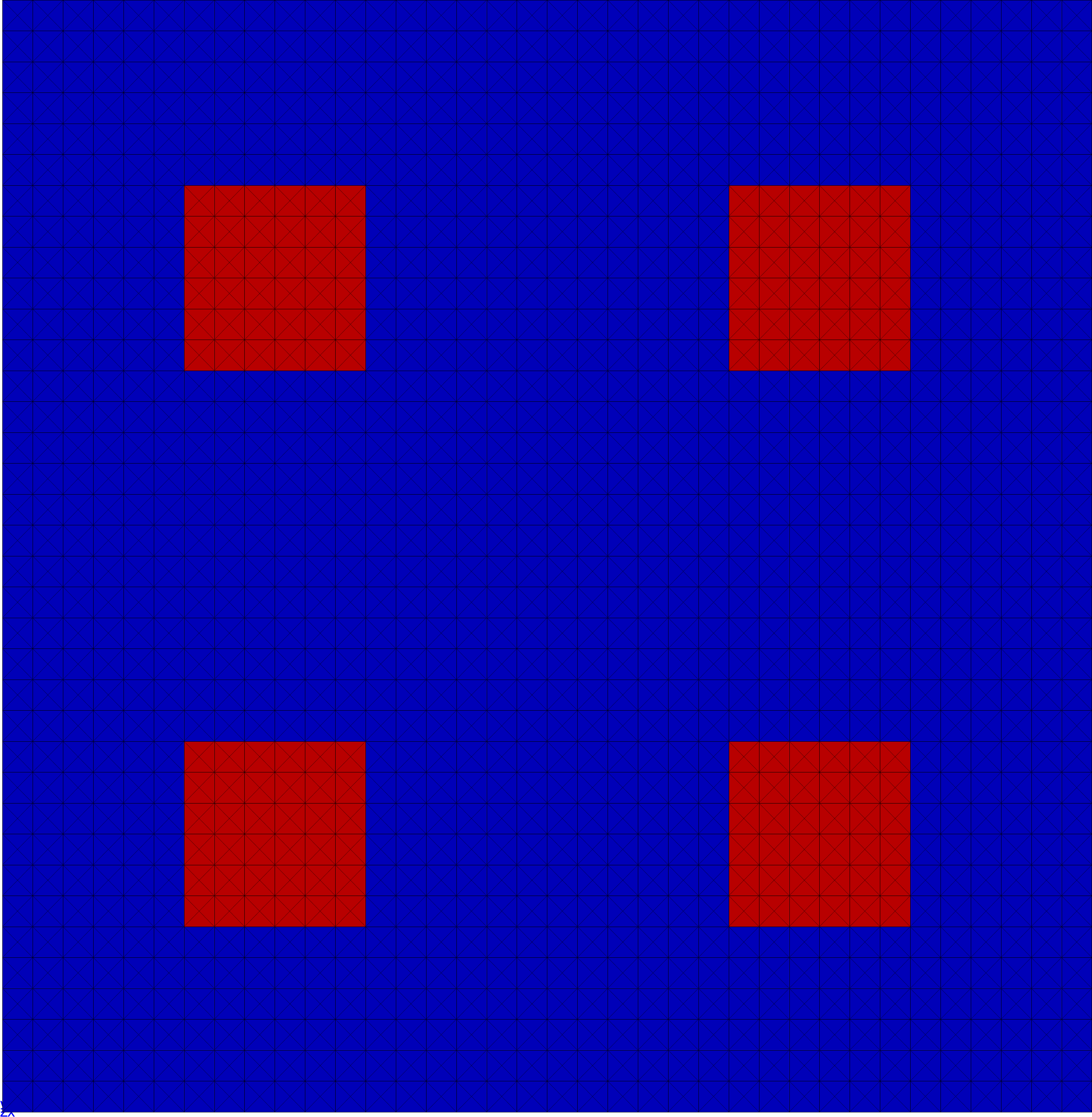}\caption{Example of mesh generated($h=L/36$)}\label{fig:heteroMesh}
\end{subfigure}\hfill
    \begin{subfigure}{.30\textwidth}
    \includegraphics[width=.18\textwidth]{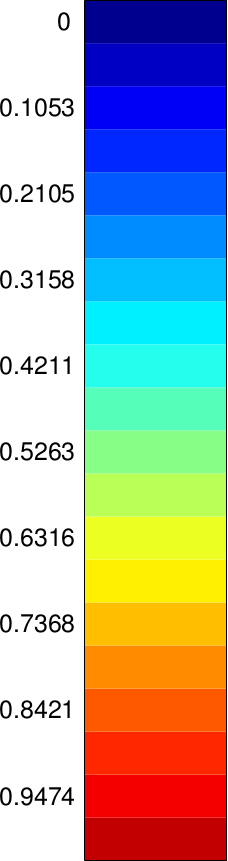} \ 
   \includegraphics[width=.78\textwidth]{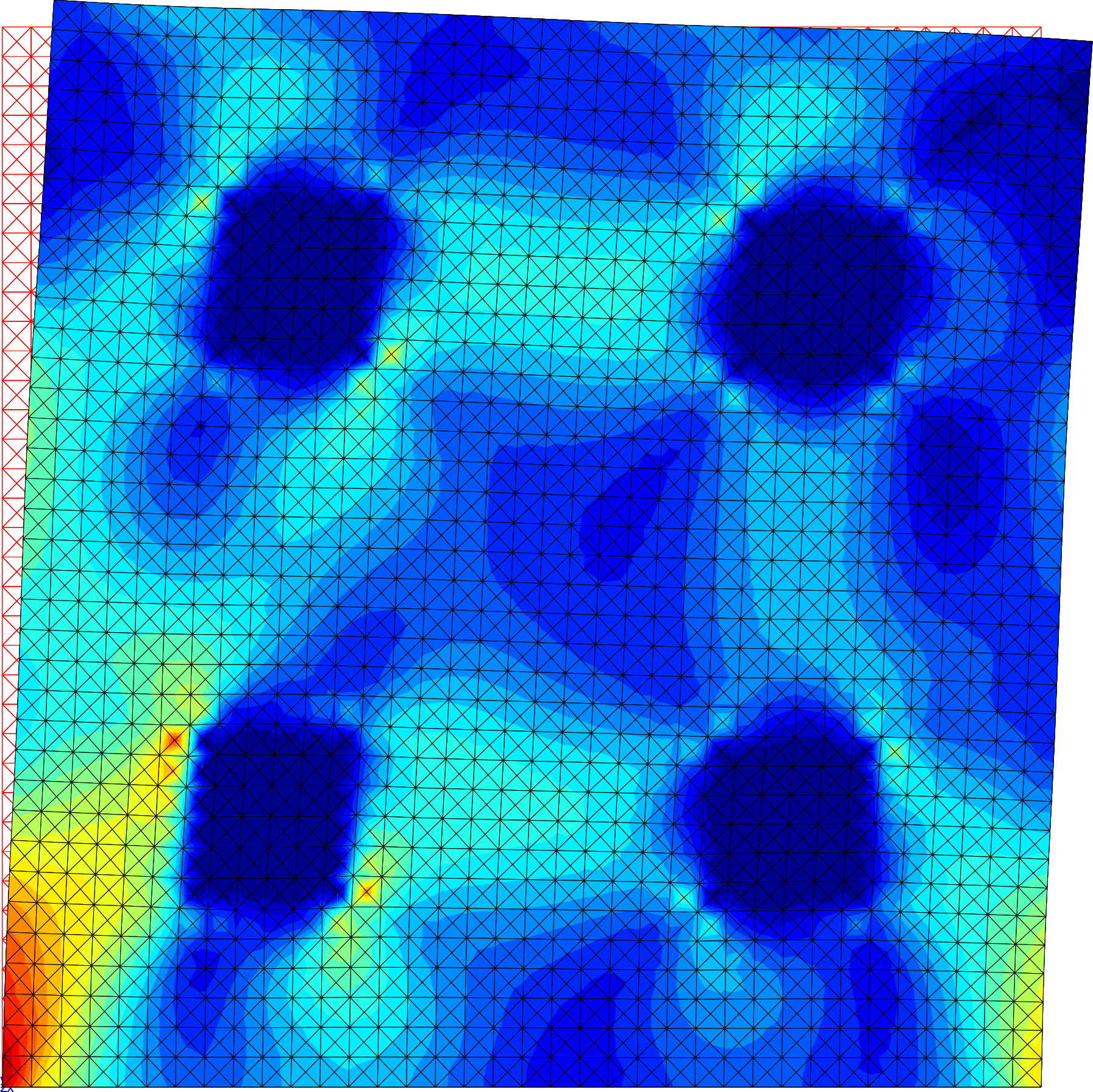}\caption{Von-Mises stresses on mesh with $h=L/36$ ($E_2=200$Pa}\label{fig:heteroMises}
    \end{subfigure}\hfill
  \caption{2D problem studied}
  \label{fig:hetero}
\end{figure}

In order to assess the performance of our parallel estimator in the presence of multiple points and heterogeneities, we consider the 2D problem of a square of side $L$ clamped on its basis, submitted to traction and shear on its upper edge while both left and right edges are assumed to be traction-free.
Four square are included in the matrix, as shown in figure \ref{fig:hetero}.
Materials are chosen to be isotropic linear elastic, with parameters $E_1=2.10^5$\,Pa and $\nu_1=\nu_2=0.3$.
$E_2$ spans values making the heterogeneity ratio varying from $10^{-6}$ to $10^6$.

A mesh constituted by P1 triangular elements of characteristic size $h=\frac{L}{36}$ is used.
A sequential computation on the whole domain is conducted, followed by domain decomposition calculations obtained by  splitting of the original domain in an increasing number $\Nsd$ of subdomains, with $\Nsd=5,9,18,36$.

Both BDD and FETI algorithms used to solve the substructured problems are respectively equipped with Neumann-Neumann 
and Dirichlet preconditioners with stiffness scaling operators. 
Beside, the convergence criterion of the solver is set to a value making the algebraic error negligible with respect to the discretization error. 
At convergence of the solver, we perform error estimation, as described in Section~\ref{sec:principles}. The method 
used to build statically admissible 
stress field is either the EET technique \citep{Lad83} or the flux-free technique \citep{Par09}, written down SPET. Element 
problems are solved  with p-refinement (three degrees higher polynomial basis \citep{Bab94}). 

\subsection{Quality of error estimators with various partitioning}
\label{sec:numassdecomp}
In this subsection, we consider four substructurings and two distributions of materials with inverse heterogeneity ratios (soft inclusions or stiff inclusions).
\subsubsection{Domain partitioning and notation}
\begin{figure}[ht]
  \centering
\begin{subfigure}{.35\textwidth}  
  \includegraphics[width=1.\textwidth]{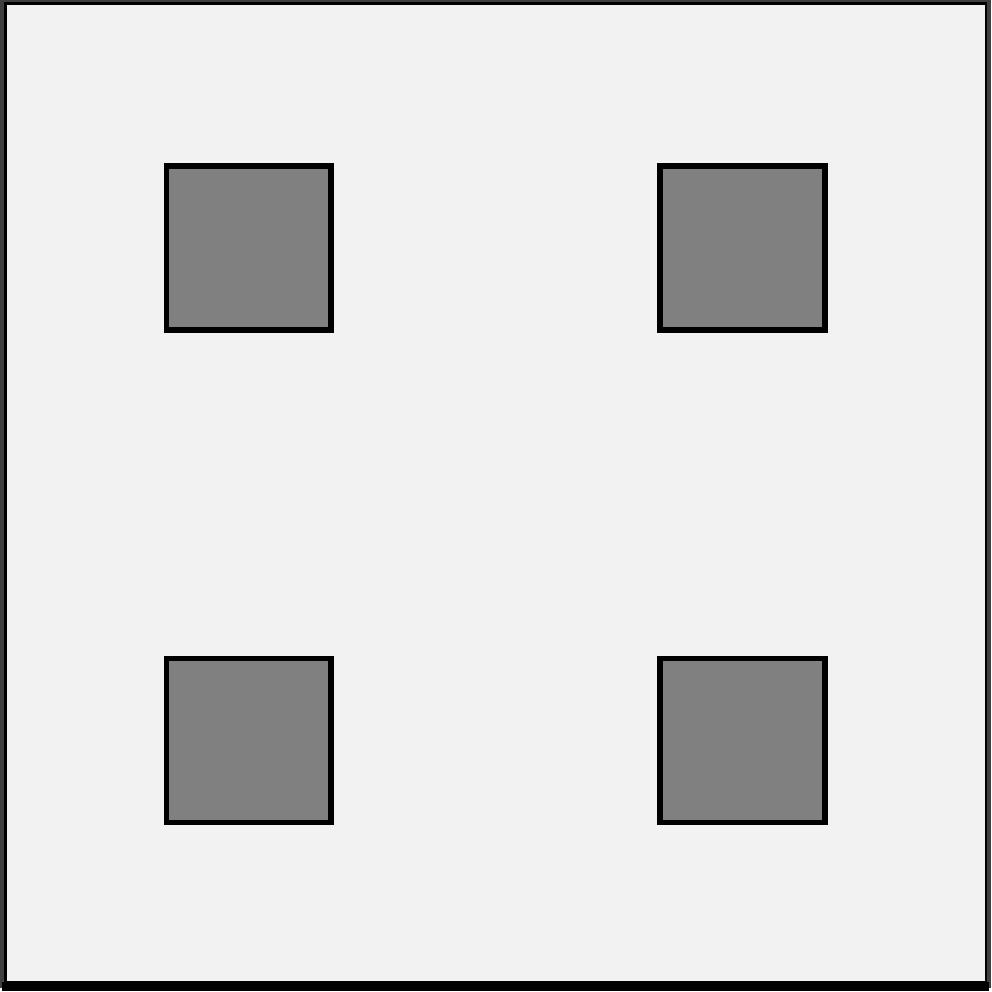}\caption{Decomposition into 5 subdomains (matrix and four inclusions)}\label{fig:dd_5}
  \end{subfigure}
\qquad
  \begin{subfigure}{.35\textwidth}  
  \includegraphics[width=1.\textwidth]{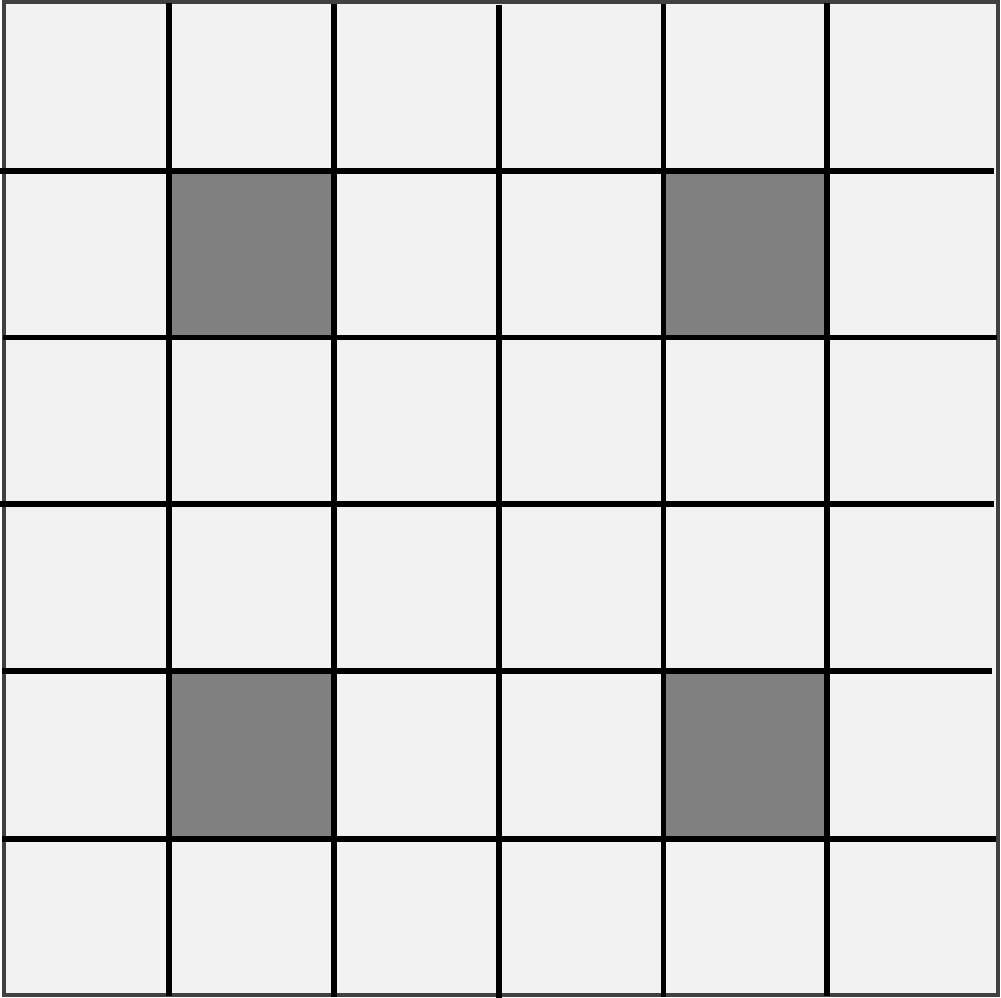}\caption{Decomposition into 36 subdomains}\label{fig:dd_36}
  \end{subfigure}
%
  \caption{Substructurings with homogeneous domains}
  \label{fig:dd_homog_decomposition}
\end{figure}

\begin{figure}[ht]
  \centering
  \begin{subfigure}{.35\textwidth}  
  \includegraphics[width=1.\textwidth]{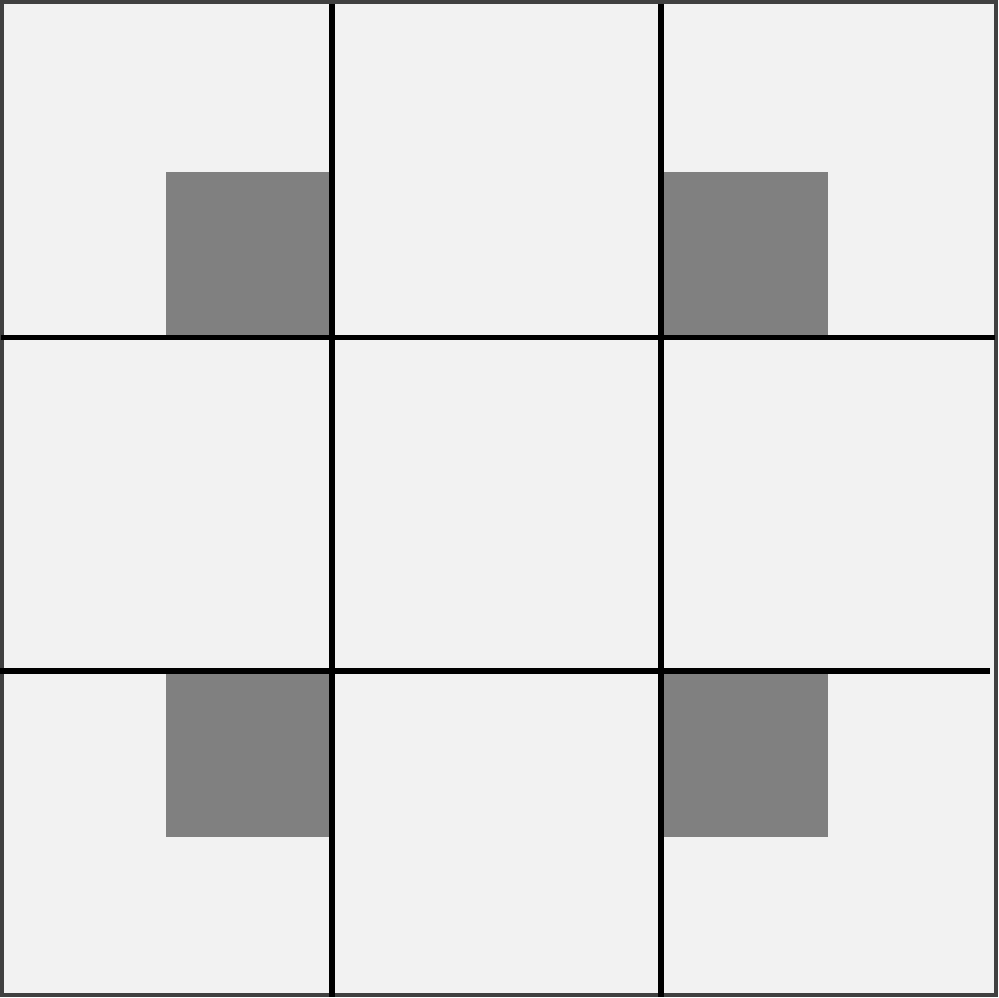}\caption{Decomposition into 9 subdomains}\label{fig:dd_9}
  \end{subfigure}
  \qquad
    \begin{subfigure}{.35\textwidth}  
  \includegraphics[width=1.\textwidth]{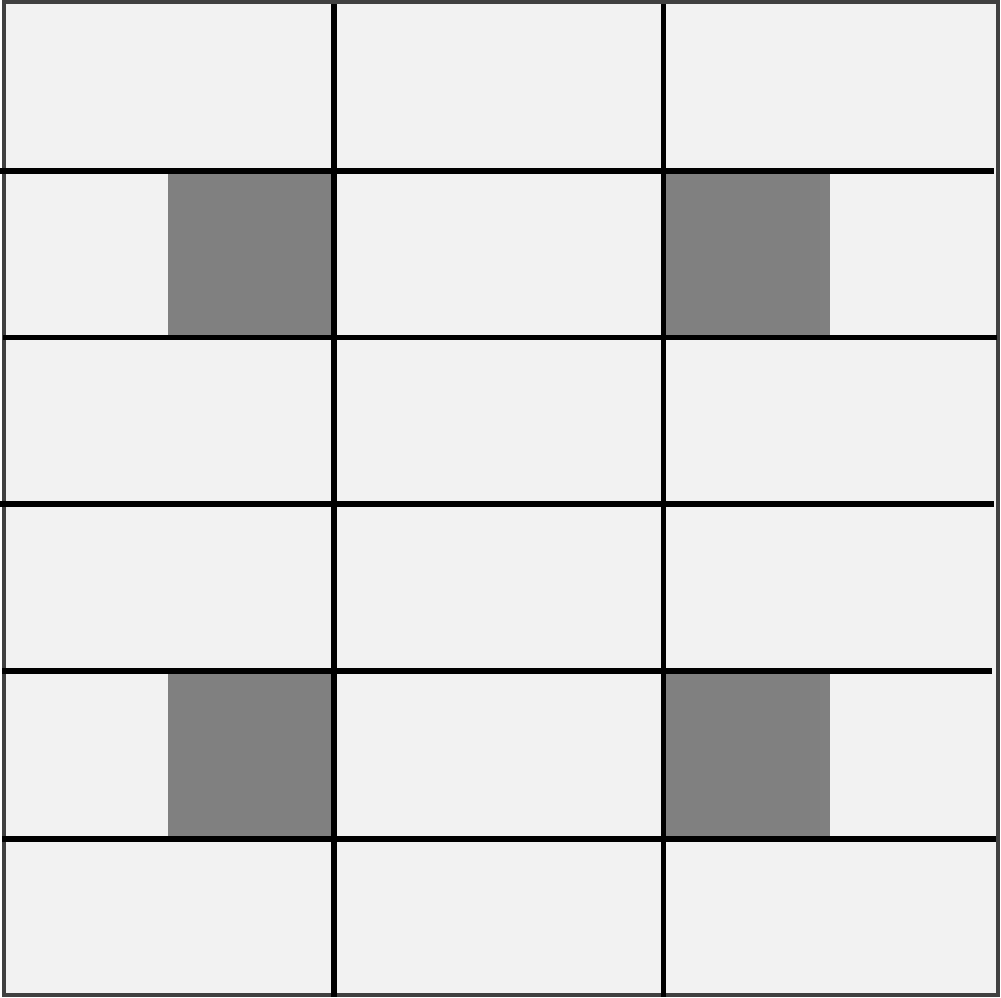}\caption{Decomposition into 18 subdomains}\label{fig:dd_18}
  \end{subfigure}
  \caption{Substructurings with heterogeneous domains}  \label{fig:dd_heter_decomposition}
\end{figure}

 In Figure~\ref{fig:dd_homog_decomposition}, we give two decompositions that involve homogeneous subdomains. The first decomposition is directly based on the position of the 
inclusions and is without multiple points. In Figure~\ref{fig:dd_heter_decomposition}, we illustrate decompositions in 9 and 18 subdomains that lead to heterogeneities inside subdomains. 
%

In the presentation of 
numerical results, subscript \textit{optim} refers to optimization with respect to multiple points in substructured 
context and optimization inside 
the EET procedure in case of heterogeneities in the domain. For instance, since the substructuring into 9 subdomains 
leads to heterogeneous subdomains, $EET_{optim}$ means that both optimizations for multiple points and for 
heterogeneity inside each subdomain have been done.

\subsubsection{Case of soft inclusions}
In Table~\ref{tab:various_dd}, we present the error estimations normalized by the energy norm 
of the finite element 
solution in the case of $\frac{E_2}{E_1}=10^{-5}$ for the four considered substructurings.

\begin{figure}[ht]
\includegraphics[scale=0.35]{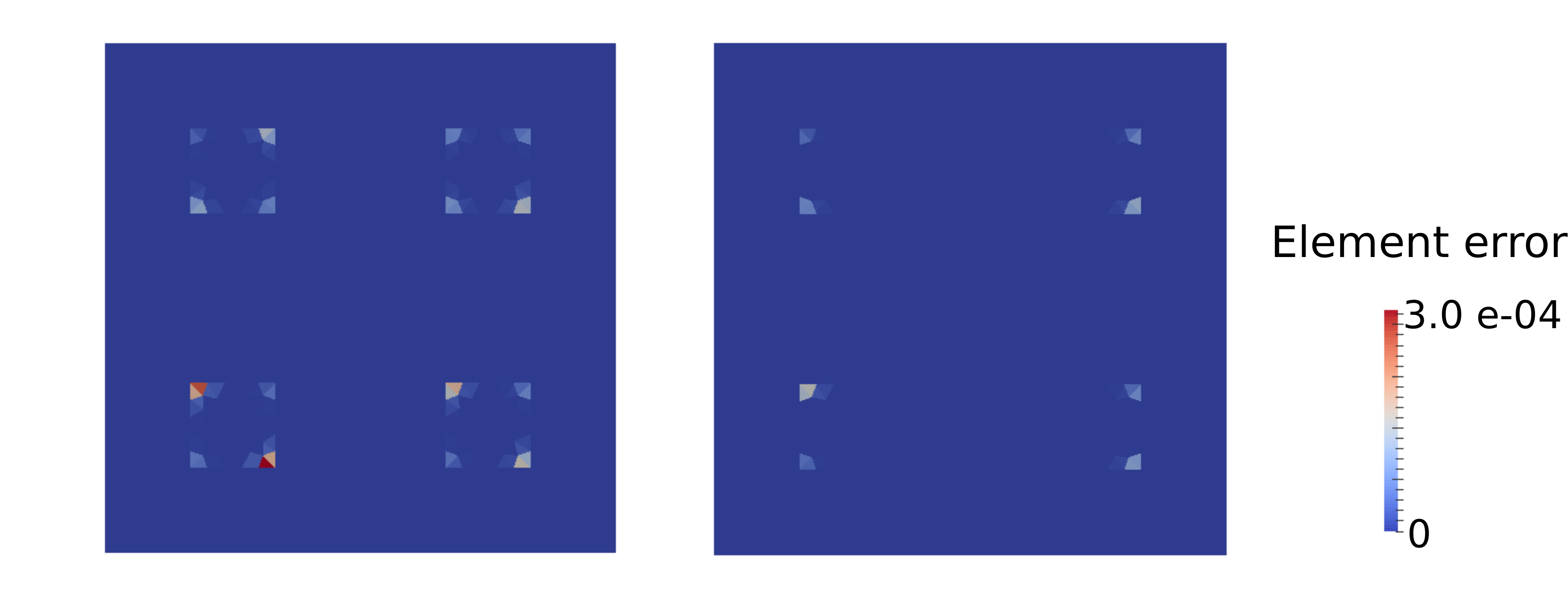}
  \caption{Error maps obtained with EET-based error estimation without optimization (left) and with optimization 
(right) with 18 subdomains}\label{fig:error_maps}
\end{figure}

\begin{table}[ht]\centering%

\begin{tabular}{|c|c|c|c|c|}
\hline 
 Number of subdomains& $EET$ &   $EET_{optim}$ &$SPET$ & $SPET_{optim}$\\ 
\hline \ 
1 & 30.5 &3.64 $10^{-1}$& 2.43 $10^{-1}$ &  \\
\hline 
5  & 3.64 $10^{-1}$&  & 2.43 $10^{-1}$& \\
\hline 
9  &  37.9 &29.6& 28.9 & 22.7\\
\hline 
18 &  40.5 &25.3& 33.7 & 19.2\\
\hline 
36  & 44.1 & 3.99 $10^{-1}$ &38.9 & 3.31 $10^{-1}$\\
\hline 
\end{tabular} 
\caption{Dependence of the error estimators wrt the substructuring in presence of soft inclusions ($10^5$ ratio of Young's moduli, values normalized by the energy).}\label{tab:various_dd}
\end{table}
 We make the following observations:
 \begin{itemize}
  \item in the sequential case, the optimization in the EET technique (presented in Section~\ref{sec:seqhetero}) 
enables to recover a good error estimation close 
to the one provided by the SPET which is always the better.
\item With 5 subdomains, there is no multiple points and domains are homogeneous. Therefore no optimization was done.
\item With 9 or 18 subdomains, the optimization at the multiple points enables to better the error estimation but does 
not lead to results as accurate as for the sequential optimization. This is due to the strong heterogeneity at the 
interface. As shown on the error maps in figure~\ref{fig:error_maps}, the optimization at the multiple points 
reduces the parasite error in those multiples points but not along the interfaces crossing heterogeneities. 
 \end{itemize}

\subsubsection{Case of stiff inclusions}
In Table~\ref{tab:various_dd2}, we present the relative error estimations for $\frac{E_2}{E_1}=10^5$ for the four 
considered 
substructurings. 

\begin{table}[ht]\centering%

\begin{tabular}{|c|c|c|c|c|}
\hline 
 Number of subdomains& $EET$ &  $ EET_{optim}$ &$SPET$ & $SPET_{optim}$\\ 
\hline \ 
1 & 2.36 $10^{-1}$& 2.23 $10^{-1}$& 1.64 $10^{-1}$& \\
\hline 
5  &  2.34 $10^{-1}$& & 1.86 $10^{-1}$&  \\
\hline 
9  & 2.43 $10^{-1}$& 2.35 $10^{-1}$& 1.81 $10^{-1}$& 1.80 $10^{-1}$ \\
\hline 
18 & 2.41 $10^{-1}$& 2.37 $10^{-1}$&1.93 $10^{-1}$ & 1.92 $10^{-1}$\\
\hline 
36  & 2.37 $10^{-1}$& 2.36 $10^{-1}$& 2.00 $10^{-1}$& 1.99 $10^{-1}$\\
\hline 
\end{tabular} 
\caption{Dependence of the error estimators wrt the substructuring in presence of stiff inclusions ($10^5$ ratio of Young's moduli, values normalized by the energy).}\label{tab:various_dd2}
\end{table}

We observe that, in the sequential case, even without optimization, the error estimation provided by the EET is 
almost as accurate as the one provided by the SPET. We notice that the optimization in the EET does not worsen the 
error estimation nor improves it a lot. The error estimation in parallel resolutions is already as accurate 
as the error 
estimation in a sequential computations. As a consequence, it is not surprising that the specific treatment of 
multiple 
point does not lead to better results.

\subsection{Quality of error estimators with increasing heterogeneity
  between subdomains}
\label{sec:numasshetero}
\begin{figure}[ht]
\centering
\begin{tikzpicture}
\begin{loglogaxis}[
scale=1.2,
xlabel=ratio $E_1/E_2$,
ylabel=relative error,
ymin=0.1,ymax=1,
legend style={at={(.05,.55)}, anchor=south west}
]
\addplot[color=blue, mark=square*] table[x=Ratio,y=EET]{dur_dans_mou_rel.txt};
\addlegendentry{\begin{scriptsize}EET\end{scriptsize}}
\addplot[color=blue, mark=square] table[x=Ratio,y=SPET]{dur_dans_mou_rel.txt};
\addlegendentry{\begin{scriptsize}SPET\end{scriptsize}}
\addplot[color=cyan, mark=square*] table[x=Ratio,y=EEToptim]{dur_dans_mou_rel.txt};
\addlegendentry{\begin{scriptsize}{EEToptim}\end{scriptsize}}
\addplot[color=red, mark=*] table[x=Ratio,y=DDEET]{dur_dans_mou_rel.txt};
\addlegendentry{\begin{scriptsize}DD EET\end{scriptsize}}
\addplot[color=red, mark=o] table[x=Ratio,y=DDSPET]{dur_dans_mou_rel.txt};
\addlegendentry{\begin{scriptsize}DD SPET\end{scriptsize}}
\addplot[color=magenta, mark=triangle*] table[x=Ratio,y=DDoptimEET]{dur_dans_mou_rel.txt};
\addlegendentry{\begin{scriptsize}DD optim EET\end{scriptsize}}
\addplot[color=magenta, mark=triangle] table[x=Ratio,y=DDoptimSPET]{dur_dans_mou_rel.txt};
\addlegendentry{\begin{scriptsize}DD optim SPET\end{scriptsize}}
\end{loglogaxis}
\end{tikzpicture}
\caption{Quality of error estimator with increasing heterogeneity : soft inclusions}
\label{fig:heterogeneity_soft}
\end{figure}
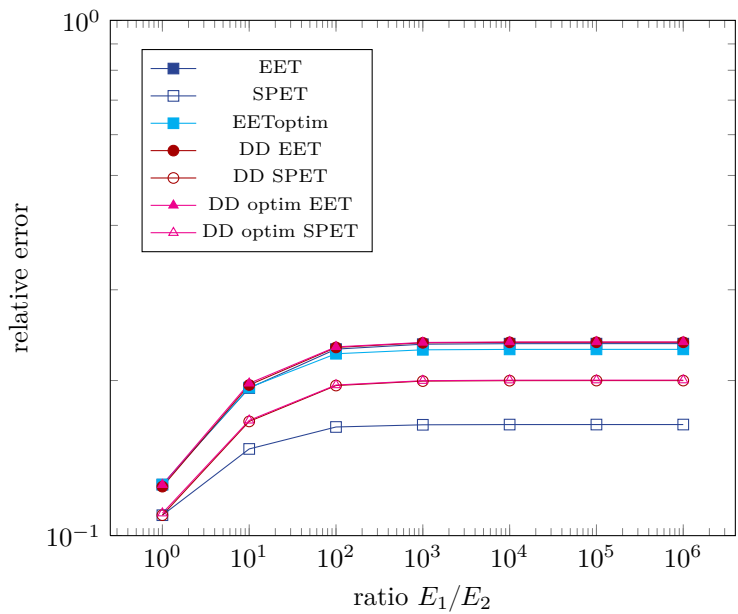

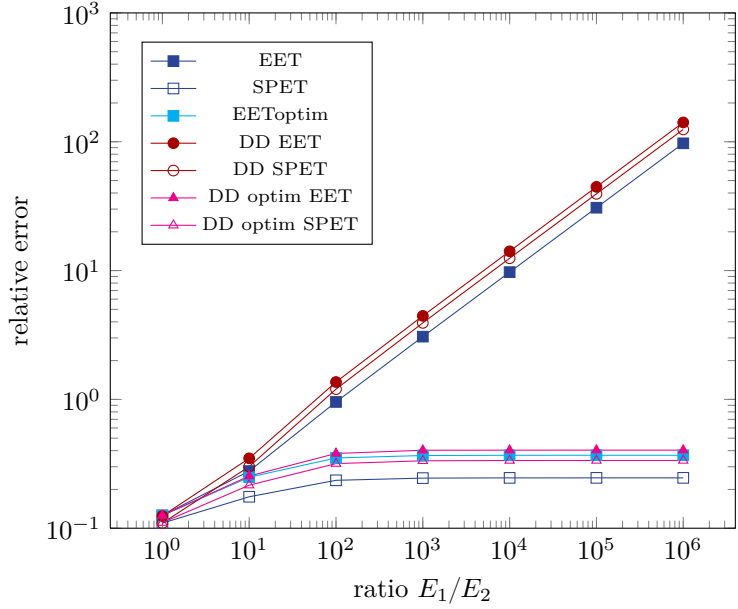
\begin{figure}
\centering
\begin{tikzpicture}
\begin{loglogaxis}[
scale=1.2,
xlabel=ratio $E_1/E_2$,
ylabel=relative error,
ymin=0.1,ymax=1e3,
legend style={at={(.05,.55)}, anchor=south west}
]
\addplot[color=blue, mark=square*] table[x=Ratio,y=EET]{mou_dans_dur_rel.txt};
\addlegendentry{\begin{scriptsize}EET\end{scriptsize}}
\addplot[color=blue, mark=square] table[x=Ratio,y=SPET]{mou_dans_dur_rel.txt};
\addlegendentry{\begin{scriptsize}SPET\end{scriptsize}}
\addplot[color=cyan, mark=square*] table[x=Ratio,y=EEToptim]{mou_dans_dur_rel.txt};
\addlegendentry{\begin{scriptsize}{EEToptim}\end{scriptsize}}
\addplot[color=red, mark=*] table[x=Ratio,y=DDEET]{mou_dans_dur_rel.txt};
\addlegendentry{\begin{scriptsize}DD EET\end{scriptsize}}
\addplot[color=red, mark=o] table[x=Ratio,y=DDSPET]{mou_dans_dur_rel.txt};
\addlegendentry{\begin{scriptsize}DD SPET\end{scriptsize}}
\addplot[color=magenta, mark=triangle*] table[x=Ratio,y=DDoptimEET]{mou_dans_dur_rel.txt};
\addlegendentry{\begin{scriptsize}DD optim EET\end{scriptsize}}
\addplot[color=magenta, mark=triangle] table[x=Ratio,y=DDoptimSPET]{mou_dans_dur_rel.txt};
\addlegendentry{\begin{scriptsize}DD optim SPET\end{scriptsize}}
\end{loglogaxis}
\end{tikzpicture}
\caption{Quality of error estimator with increasing heterogeneity : stiff inclusions}
\label{fig:heterogeneity_stiff}
\end{figure}

In this subsection, we decompose the structure into 36 identical square subdomains (see 
Figure~\ref{fig:dd_36}) and we keep the same 
substructuring and discretization. Note that with this decomposition, every subdomain is homogeneous. We 
study the 
behavior of the error estimators with increasing heterogeneity.
 In the presentation of the results, we adopt the 
following notations : 

\begin{itemize}
 \item EET : sequential resolution and use of the EET procedure for the construction of admissible stress fields
 \item SPET : sequential resolution and use of the SPET procedure for the construction of admissible stress fields
 \item EET optim : sequential resolution and use of the EET procedure with optimization for the construction of 
admissible stress fields
  \item DD EET : parallel resolution  and use of the EET procedure
 \item DD SPET : parallel resolution  and use of the SPET procedure
 \item DD optim EET : parallel resolution with optimization on the multiples points and use of the EET procedure
  \item DD optim SPET : parallel resolution with optimization on the multiples points and use of the SPET procedure
\end{itemize}
In Figure~\ref{fig:heterogeneity_soft}, $E_2$ is smaller than $E_1$ and in Figure~\ref{fig:heterogeneity_stiff}, $E_2$ 
is 
larger than $E_1$.
 
We observe that in the cases where basic estimators behave poorly (EET, DD EET and DD SPET with soft inclusions), the 
optimization enables to recover correct order of magnitude whatever the heterogeneity ratio. In other cases all estimators give quite close results.


\section{Application to error estimation in FETI-DP algorithm}
\label{sec:fetidp}
The results of previous sections are very useful in order to derive 
a procedure to obtain admissible fields in the FETI-DP method \cite{Far00bis, Far01} and to prove that the bounds 
with separated contributions of \cite{VREY.2013.1, Rey15, Rey16} also apply to this algorithm.

Note that all results presented here are also valid for BDDC \cite{Doh03} which corresponds to the derivation of the same ideas as FETI-DP starting from the primal approach.

\subsection{The FETI-DP algorithm}
In the FETI-DP algorithm, primal constraints are incorporated in the formulation so that subdomains always remain 
weakly connected. The chosen connections make all local problems well posed, and  act like a coarse problem which 
warranties the scalability of the method even in difficult situations: \cite{Doh03, Kla07, 
Klawonn:2008:AFA:1405063.1405075}. The most basic constraint consists, in 2D, in ensuring the continuity of the 
displacement at corners (which are the boundary of edges, multiple points being a particular class of such nodes). This 
is the type of constraints we will focus on subsequently.

Note that more general constraints can be chosen, like weighted averages over faces (in 3D).
Often these constraints are 
implemented using local change of basis (so that even complicated constraints end up to be applied on single modified 
degrees of freedom). Anyhow these changes of basis do not modify the methodology to compute admissible fields because 
continuity and balance are properties independent from the chosen basis. 

Let the subscript $c$ represent corners, and the subscript $r$ represent the remaining degrees of freedom of which we 
distinguish the internal dofs $i$ from the other interface dofs $o$ (not corners). We define in a straightforward 
manner the restriction of trace and assembly operators on corners and other degrees of freedom. The FETI-DP problem can 
be written as: 
\begin{equation}\text{Find }\Lam_o\text{ such that }
\left\{ \begin{aligned}
&\left\{ 
\begin{aligned}&\stiff\s\dep\s=\force\s + \traceh_o\sT\dassem_o\sT \Lam_o \text{, }\forall s \\
& \text{under the constaint } \sum_s \dassem\s_c \traceh_c\s\dep\s =0 \end{aligned} \right.\\
&\text{such that } \sum_s\dassem_o\s\traceh_o\s\dep\s=0
\end{aligned} \right.
\end{equation}
Algorithm~\ref{alg:fetidp} presents the associated algorithm, with the following notations:
\begin{itemize}
\item  $\stiff_{cc}^*$ matrix of the coarse problem:
\begin{equation}
 \stiff_{cc}^*= \sum_s\passem_c\s\left( \stiff_{cc}\s  - 
 \stiff_{rc}\sT  \stiff_{rr}\sinv \stiff_{rc}\s \right)
\passem_c\sT
\end{equation}
\item Scaled assembly matrix $\widetilde{\dassem}_o\s$, such that $\sum\dassem_o\s\widetilde{\dassem}_o\sT=\mathbf{I}$
\item Forward problem $(\dep\s,\lam_c\s)=\text{Solve}_L(\lam_o\s,\force\s)$:
\begin{equation}\label{eq:fetidpforward}
\left\{ \begin{aligned}
&\dep\s=  \begin{pmatrix}\dep\s_{1,r} + \dep\s_{3,r} \\ \dep_c\s
\end{pmatrix},\quad \lam_c\s=\stiff_{cr}\s (\dep\s_{1,r} + \dep\s_{3,r}) + \stiff\s_{cc} \dep\s_c \\
&\text{ with }  \dep\s_{1,r} \text{ solution to } \stiff\s_{rr} \dep\s_{1,r}=\traceh\sT\lam_o\s+\force\s\\
& \text{ and } \dep\s_c =  -\passem\sT_c {\stiff_{cc}^*}^{-1}
\sum_j \passem^{(j)}_c \stiff^{(j)}_{cr} \dep^{(j)}_{1,r} \\
&\text{ and } \dep\s_{3,r} \text{ solution to } \stiff\s_{rr} \dep\s_{3,r}=-\stiff\s_{rc} \dep\s_c  \end{aligned}
\right.
\end{equation}
\item Preconditioning $(\tlam\s,\delta\dep\s)=\text{Solve}_S(\dep_{o}\s)$
\begin{equation*}
\left\{ \begin{aligned}
&\stiff\s \delta\dep\s =  \eft\sT \tlam\s \\
&\text{under the conditions } \eft\s_{o} \delta\dep\s=\dep_{o}\s,\  \eft\s_{c} \delta\dep\s=0
\end{aligned}
\right.
\end{equation*}
\end{itemize}

\begin{algorithm2e}[ht]\caption{FETI-DP: main unknown $\Lam_o$}\label{alg:fetidp}
Initialization $\Lam_o$ \;  %
$(\dep_N\s,\lam_c\s)=\text{Solve}_L(\dassem\sT_o \Lam_o,\force\s)$\tcp*[r]{$\lam\s_N = \left[\dassem\sT_o \Lam_o;\lam_c\s\right]$}
Compute residual $\res=( \sum_s\dassem_o\s \traceh_o\s\dep_N\s )$\;
Define local displacement $\delta\dep\s_{o}=\widetilde{\dassem}_o\sT
\res$ \;
$(\tlam\s,\delta\dep\s)=\text{Solve}_{S}(\delta\dep_{o}\s)$  
\tcp*[r]{$\begin{aligned}&\dep_D\s=\dep_N\s-{\delta\dep}\s \\
&\lam_D\s=\lam_N\s-\tlam\s\end{aligned}$}
Compute preconditioned residual $\bz=\sum_s\widetilde{\dassem}_o\s\tlam_o\s$ \;
Compute search direction $\bw=\bz$\;
\While{$\sqrt{\res^T\bz}>\epsilon$}{%
  $(\delta\dep_N\s,\Delta\lam_c\s)=\text{Solve}_L(\dassem\sT_o   \bw,0)$\;
  $\bp=\sum_s\dassem_o\s \traceh\s_o\delta\dep_N\s$\;
  $\alpha=(\res^T\bz)/(\bp^T\bw)$\;
  $\Lam_o\leftarrow \Lam_o+\alpha \bw$  
\tcp*[r]{$\begin{aligned}&\dep_N\s\leftarrow \dep_N\s+\alpha 
\delta\dep_N\s\\&\lam_N\s=\left[\dassem\sT_o \Lam_o; \lam\s_c+\alpha\Delta\lam\s_c\right]\end{aligned}$}
  $\res \leftarrow \res-\alpha \bp$\;
  $\delta\dep_{o}\s=\widetilde{\dassem}_o\sT \res$\;
  $(\tlam\s,\delta\dep\s)=\text{Solve}_{S}(\delta\dep_{o}\s)$ 
\tcp*[r]{$\begin{aligned}&\dep_D\s=\dep_N\s-{\delta\dep}\s \\
&\lam_D\s=\lam_N\s-\tlam\s\end{aligned}$}
  $\bz = \sum_s\widetilde{\dassem}_o\s\tlam_o\s$\;
  $\bw \leftarrow \bz - (\bp^T\bz)/(\bp^T\bw) \bw $
}%
\end{algorithm2e}

\subsection{Parallel reconstruction of admissible fields and error estimation}

In algorithm~\ref{alg:fetidp}, we showed how it was possible at no cost to build the fields 
$\dep_N\s,\dep_D\s,\lam_N\s$. As proved in the following paragraphs, the following properties hold:
\begin{itemize}
\item $\dep_D\s$ defines  a continuous field of displacement:
\begin{equation}
\sum_s\dassem_o\s\traceh_o\s {\dep_{D}\s} =  \sum_s\dassem_o\s \traceh_o\s \dep_{N}\s -  \sum_s\dassem_o\s\traceh_o\s
\widetilde{\dassem}_o\sT   \dassem_o\s\traceh_o\s \dep_{N}\s
= 0 
\end{equation}
The continuity at corners is at the basis of the coarse problem and is not perturbed during the preconditioning (where null displacement is imposed on corners)
\item $(\dep_N\s,\lam_N\s)$ satisfy the subdomains' equilibrium:
\begin{equation}
\stiff\s\dep_N\s=\force\s+\traceh\sT\lam_N\s
\end{equation}
and reactions are balanced at the interface:
\begin{equation}
\begin{aligned}
&\sum_s\passem_o\s\lam\s_{N,o}= \sum_s\passem_o\s\dassem_o\s \Lam_o = 0 \\
&\sum_s\passem_c\s\lam\s_{N,c}= 0 \text{ from definition } \eqref{eq:fetidpforward}  \\
\end{aligned}
\end{equation}
\end{itemize}
This means that $u_D=(\shapev\s\dep_D\s)$ is a kinematically admissible displacement field, and that $\sig_N\s=\hooke:\eps(\shapev\s\dep_N\s)$ can be used, together with $\lam_N\s$, to compute a globally statically admissible stress field $\hsig_h$.

Once those fields are 
obtained, the error estimation consists in the evaluation of the quantity 
\begin{equation*}
  \left.\ecr{\dom}(\hu_h,\hsig_h)\right. = \sqrt{\sum_{s=1}^{\Nsd} 
\left(\ecr{\dom\s}(\hu\s_h,\hsig\s_h)\right)^2.}
\end{equation*}
which is a strict upper bound of the true error.

\subsection{Separation of the sources of the error}
We have the following property, using notations of algorithm~\ref{alg:fetidp}:
\begin{equation}
\begin{aligned}
\vvvert {u}_N-{u}_D \vvvert_{\Omega}^2 &= \sum_s \int_{\Omega\s} \eps(\un{u}_D\s-\un{u}_N\s) : \hooke: \eps(\un{u}_D\s-\un{u}_N\s) d\Omega\\
&= 
 \sum_s \left(\lam_D\s-\lam_N\s\right)^T{\traceh\s} \left(\dep_D\s-\dep_N\s\right)\\
 &=\sum_s\tlam\sT {\traceh\s} \tdep\s = \sum_s \tlam_o\sT\tdep_o\s \\& =  \mathbf{r}^T \mathbf{z}  \\ 
\end{aligned}
\end{equation}
where we have used the fact that $\tdep_c\s=0$

This means that the quantity $\vvvert {u}_N-{u}_D \vvvert_{\Omega}^2$ is naturally computed at each iteration. Thus all 
the results involving separation of sources of error \eqref{eq:sepsource} can be used in FETI-DP (see also  
\cite{Rey15, Rey16} for bounds on quantities of interest and lower bounds).

\section{Conclusion}

In this article, we emphasize the role of multiple points for the error estimation in the framework of domain decomposition methods. We propose an optimization to reconstruct nodal reactions at these specific points that is optional in absence of heterogeneity but is necessary to recover a efficient error estimator in case of strong heterogeneities. 
By the way, because of strong similarities between optimization at multiple points and reconstruction of fluxes in 
star-patches in the EET algorithm, we propose an enhancement of the EET in case of strong heterogeneities. The 
numerical assessments on a 2D mechanical structure show the quality of 
the proposed improvements. Finally, the specific attention towards the multiple points also allows us to extend the 
parallel error estimation 
procedure to FETI-DP (and BDDC) algorithms.

\subsection*{Acknowledgement}
The authors wish to thank professor Nicolas Mo\"es for pointing out the particular role of multiple points in the recovery process.

\bibliography{Biblio}







\end{document}